\documentclass[aop,MSNbibl,seceqn,dvips]{arximspdf}
\usepackage{graphicx}
\doi{10.1214/12-AOP794} 
\volume{42}
\issue{3}
\pubyear{2014}
\firstpage{1020}
\lastpage{1053}

\makeatletter
\def\l{}
\newcommand{\rrVert}{\Vert}
\newcommand{\rrvert}{\vert}
\newcommand{\llVert}{\Vert}
\newcommand{\llvert}{\vert}
\def\cal{\mathcal}

\newtheorem{Lemma}{Lemma}[section]
\newtheorem{proposition}[Lemma]{Proposition}
\newtheorem{theorem}[Lemma]{Theorem}
\newproclaim{remark}[Lemma]{Remark}
\newproclaim{example}[Lemma]{Example}

\newcommand{\reals}{\mathbb{R}}
\newcommand{\bbr}{\reals}
\newcommand{\vep}{\varepsilon}

\def\definedas{\stackrel{\Delta}{=}}
\newcommand{\BX}{\mathbf{X}}
\newcommand{\bt}{\mathbf{t}}
\newcommand{\bs}{\mathbf{s}}
\newcommand{\ba}{\mathbf{a}}
\newcommand{\bb}{\mathbf{b}}
\newcommand{\bnull}{\mathbf{0}}
\newcommand{\bomega}{\bolds{\omega}}
\newcommand{\one}{\mathbf{1}}

\newcommand{\cov}{\operatorname{cov}}

\newcommand{\cH}{\mathcal{H}}
\newcommand{\cP}{\mathcal{P}}
\newcommand{\CXa}{\mathcal{C}_\BX(\ba)}
\newcommand{\CXab}{\mathcal{C}_\BX(\ba,\bb)}
\newcommand{\CXxi}{\mathcal{C}_\BX(\ba; \xi)}
\newcommand{\CXbxi}{\mathcal{C}_\BX(\ba,\bb; \xi)}
\newcommand{\Wxi}{\mathcal{W}_\xi}
\newcommand{\Wa}{\mathcal{W}_a}
\newcommand{\Wab}{\mathcal{W}_{a,b}}

\def\definedby{\stackrel{\Delta}{=}}
\newcommand{\eqref}[1]{(\ref{#1})}
\makeatother

\begin{document}
\begin{frontmatter}

\runtitle{Paths in Gaussian excursions}
\title{On the existence of paths between points in high level
excursion sets of Gaussian random fields}

\begin{aug}
\author[A]{\fnms{Robert J.}~\snm{Adler}\thanksref{t1,t4}\ead[label=e1]{robert@ee.technion.ac.il}\ead[label=u1,url]{webee.technion.ac.il/people/adler}},
\author[A]{\fnms{Elina}~\snm{Moldavskaya}\thanksref{t3,t5}\ead[label=e3]{elinamoldavskaya@gmail.com}}
\and\break
\author[B]{\fnms{Gennady}~\snm{Samorodnitsky}\corref{}\thanksref{t1,t2}\ead[label=e2]{gs18@cornell.edu}\ead[label=u2,url]{www.orie.cornell.edu/\textasciitilde gennady/}}
\runauthor{R.~J. Adler, E.~Moldavskaya and G.~Samorodnitsky}
\affiliation{Technion, Technion and Cornell University}
\thankstext{t1}{Supported in part by US-Israel Binational
Science Foundation, 2008262.}
\thankstext{t2}{Supported in part by ARO Grant W911NF-07-1-0078, NSF
Grant DMS-1005903 and NSA Grant
H98230-11-1-0154 at Cornell University.}
\thankstext{t3}{Supported in part by Israel Science Foundation, 853/10.}
\thankstext{t4}{Supported in part by AFOSR FA8655-11-1-3039.}
\thankstext{t5}{Supported in part by Office for Absorption of
New Scientists and
EC-FP7-230804 Grant.}
\address[A]{R.~J. Adler\\
E. Moldavskaya\\
Electrical Engineering\\
Technion, Haifa\\
Israel 32000\\
\printead{e1}\\
\phantom{E-mail:\ }\printead*{e3}\\
\printead{u1}}
\address[B]{G. Samorodnitsky\\
ORIE\\
Cornell University\\
Ithaca, New York 14853\\
USA\\
\printead{e2}\\
\printead{u2}}
\end{aug}

\received{\smonth{3} \syear{2012}}
\revised{\smonth{8} \syear{2012}}


\begin{abstract}
The structure of Gaussian random fields over high levels is a well
researched and
well understood area, particularly if the field is smooth.
However, the question as to whether or not two or more points
which lie in an
excursion set belong to the same connected component has
constantly eluded analysis. We study this problem from the
point of view of large deviations, finding the asymptotic
probabilities that two such points are connected by a path
laying within the excursion set,
and so belong to the same component. In addition, we
obtain a characterization and descriptions of the
most likely paths, given that one exists.
\end{abstract}

%
\begin{keyword}[class=AMS]
\kwd[Primary ]{60G15}
\kwd{60F10}
\kwd[; secondary ]{60G60}
\kwd{60G70}
\kwd{60G17}
\end{keyword}
\begin{keyword}
\kwd{Gaussian process}
\kwd{excursion set}
\kwd{large deviations}
\kwd{exceedence probabilities}
\kwd{connected component}
\kwd{optimal path}
\kwd{energy of measures}
\end{keyword}

\end{frontmatter}

\section{Introduction}\label{secIntro}

Let $\BX= (X(\bt), \bt\in\bbr^d)$ be a real-valued
sample continuous Gaussian random field. Given a level $u$, the
excursion set of $\BX$ above the level $u$ is the random set
%
\begin{equation}
\label{eexcset} A_u= \bigl\{ \bt\in\bbr^d\dvtx X(\bt)>u
\bigr\} .
\end{equation}
Understanding the structure of the excursion sets of random fields is
a mathematical problem with many applications, and it has generated
significant interest, with several recent books on the subject (e.g.,
\cite{adlertaylor2007} and \cite{azaiswschebor2009}) and with
considerable emphasis on the topology of these
sets. One very natural question in this setting which has until now
eluded solution but which we
study in this
paper is the following: given
that two points in $\bbr^d$ belong to the excursion set, what is the
probability that they belong to the same path-connected component of the
excursion set? Specifically, let $\ba, \bb\in\bbr^d$,
$\ba\not=\bb$. Recall that a path in $\bbr^d$ connecting $\ba$ and
$\bb$ is a continuous map $\xi\dvtx  [0,1]\to\bbr^d$ with $\xi(0)=\ba$,
$\xi(1)=\bb$. We denote the collection of all such paths by $\cal P
(\ba,\bb)$ and are interested in the conditional probability
\[
\label{econdprob}
P \bigl( \exists\xi\in{\cal P}(\ba,\bb)\dvtx X \bigl(\xi(v) \bigr)>u,
\mbox{ for all } 0\leq v\leq1 | X(\ba)>u, X(\bb)>u \bigr).
\]
It is straightforward to check that we are considering measurable
collections of outcomes, so this probability is well defined.

Of course, the conditional probability
above is a
ratio of two probabilities, the denominator being no more than a
bivariate Gaussian
probability, which is well
understood. Therefore, we will concentrate on the
unconditional probability
%
\begin{equation}
\label{epathprob} \Psi_{\ba,\bb}(u) \definedas P \bigl( \exists\xi\in{\cal
P}(\ba,\bb)\dvtx X \bigl(\xi(v) \bigr)>u, \mbox{ for all } 0\leq v\leq1 \bigr).
\end{equation}
If the random field is stationary, we may, without loss of
generality, assume that $\bb=\bnull$, in which
case we will use the simpler notation $\Psi_\ba$ in
\eqref{epathprob}.

When the domain of a random field is restricted to a (compact) subset
$T\subset\bbr^d$, the points $\ba$ and $\bb$ will be assumed to be in
$T$, and the entire path in \eqref{epathprob} will be
required to lie in $T$ as well (the implicit assumption being that $T$
contains \textit{some} path between $\ba$ and $\bb$).
Nevertheless, we will use the same notation and also write
\begin{eqnarray*}
\Psi_{\ba,\bb}(u) = P \bigl(\exists\xi\in{\cal P}(\ba,\bb) \dvtx \xi(v)\in
T \mbox{ and } X \bigl(\xi(v) \bigr)>u, \mbox{ for all } 0\leq v\leq1 \bigr).
\end{eqnarray*}
Which of the two interpretations of $\Psi_{\ba,\bb}$ is intended at
any point will be clear from the context.

We will study the logarithmic behavior of
the probability $\Psi_{\ba,\bb}(u)$ for high levels~$u$, that is, as
$u\to\infty$. We start with a large deviations approach, which, as
usual, will not only
describe the probability but also give us insight into the highest
probability configurations.
This makes up Sections \ref{secLD} and \ref{secsmallenergy}, which
follow a brief
technical Section \ref{secbackground} collecting some results on the
reproducing kernel
Hilbert space of a Gaussian process. Throughout we will treat the
general and the stationary
cases in parallel, but separately, since the stationary case is
somewhat more transparent
and more readily provides illustrative and illuminating special cases.
In particular, we will
look at a number of one-dimensional examples in Sections~\ref
{seconedim}--\ref{longsec}, where we can
compute quite a lot. Even in this case the results are new and rather
unexpected. We look
at the multidimensional case in Section \ref{secmultcase}. While
this section also contains
some interesting and surprising examples, it turns out that typical
examples involve nonconvex optimization
problems that we do not, at this stage, know how to solve in general.

\section{Some technical preliminaries}
\label{secbackground}

In this section we introduce much of the notation we will use in the
rest of the paper and recall certain important notions, concentrating
in particular on
the reproducing kernel Hilbert (RKHS) space of a Gaussian process.

Our main reference for the RKHS is van der Vaart and van Zanten
\cite{vandervaartvanzanten2008}, and we use it selectively so as to prepare
the background for using the large deviations theory of
Deuschel and Stroock \cite{deuschelstroock1989}.
An alternative route would be to have followed the new notes by Lifshits
\cite{lifshits2012}.

We consider a real-valued centered continuous Gaussian
random field $\BX= (X(\bt), \bt\in\bbr^d)$. When needed
(particularly, in the nonstationary case) we may restrict the domain
of the random field to a compact subset $T$ of $\bbr^d$. We denote
the covariance function of $\BX$ by
$R_\BX(\bs,\bt)=
\cov(X(\bs),X(\bt))$.

As is customary, when the random field is stationary, we will use the
single variable notation $R_\BX(\bt)=R_\BX(\bnull,\bt)$ for the
covariance function. In this case we denote the spectral measure
of $\BX$ by $F_\BX$, this being the symmetric,
finite, Borel probability measure on $\bbr^d$ satisfying
%
\begin{equation}
\label{espectral} R_\BX(\bt) = \int_{\bbr^d}
e^{i(\bt,\mathbf{x})} F_\BX(d\mathbf {x}),\qquad \bt\in\bbr^d .
\end{equation}

If $\BX$ is stationary, then this and local
boundedness imply that
\[
\lim_{\|\bt\|\to\infty}\frac{X(\bt)}{\|\bt\|} = 0
\]
with probability 1, so that almost all the sample paths of $\BX$
belong to the space
\[
C_0 = \Bigl\{ \bomega= \bigl(\omega(\bt), \bt\in\bbr^d
\bigr), \mbox { continuous, such that } \lim_{\|\bt\|\to\infty}\omega(\bt)/\|\bt \| = 0
\Bigr\}.
\]
Equipped with the norm
%
\begin{equation}
\label{215} \|\bomega\|_{C_0} = \sup_{\bt\in\bbr^d} \frac{|\omega(\bt
)|}{1+\|\bt\|} ,
\end{equation}
$C_0$ becomes a separable Banach space, with dual space
\[
C_0^* = \biggl\{ \mbox{finite signed Borel measures $\mu$ on $
\bbr^d$ with } \int_{\bbr^d} \|\bt\| \|\mu\|(d\bt)<
\infty \biggr\}.
\]
We view the stationary random field
$\BX$ as a Gaussian random element of $C_0$, generating a Gaussian probability
measure $\mu_\BX$ on that space.

In the absence of stationarity, we will usually consider a continuous Gaussian
random field $\BX= (X(\bt), \bt\in T)$, for a compact set
$T\subset\bbr^d$. In that case we view the random field
$\BX$ as a Gaussian random element in the space $C(T)$ of continuous functions
on $T$, equipped with the supremum norm, thus generating a Gaussian
probability measure $\mu_\BX$ on $C(T)$.

The reproducing kernel Hilbert space (henceforth RKHS)
$\cH$ of the Gaussian measure $\mu_\BX$ (or of the random field $\BX$)
is a subspace of $C_0$ or $C(T)$, depending on the parameter space of
$\BX$,
obtained as
follows. In the general case we identify~$\cH$ with the closure
$\mathcal L$
in the mean square norm of
the space of finite linear combinations $\sum_{j=1}^k a_jX(\bt_j)$ of
the values of the process, $a_j\in\bbr, \bt_j\in\bbr^d$ (or
$T$) for $j=1,\ldots, k$, $k=1,2,\ldots$
via the injection ${\mathcal L}\to C(T)$ given by
%
\begin{equation}
\label{eembedGen} H\to w_H = \bigl( E \bigl( X(\bt)H \bigr), \bt\in
T \bigr) .
\end{equation}

When $\BX$ is stationary, the RKHS $\cH$ can also be identified
with the subspace of functions, with even real parts and odd
imaginary parts, of the $L^2$ space of the spectral measure $F_\BX$
in \eqref{espectral},
via the injection $L^2(F_\BX)\to C_0$ given by
%
\begin{equation}
\label{eembedd} h\to S(h)= \biggl( \int_{\bbr^d} e^{i(\bt,\mathbf{x})}
\bar h(\mathbf{x}) F(d\mathbf{x}), \bt\in\bbr^d \biggr) .
\end{equation}

We denote by $(\cdot,\cdot)_\cH$ and $\| \cdot\|_\cH$ the inner
product and the norm in the RKHS~$\cH$. Since both injections described
above are
isometric, we have the important equalities:
%
\begin{equation}
\label{normequalities} E \bigl(H^2 \bigr) = \|w_H
\|^2_\cH.
\end{equation}
In the stationary case, these can be written somewhat more
informatively as
%
\begin{eqnarray}
\label{statnormequalities} \|h\|_{L^2(F_\BX)}^2 = \int
_{\bbr^d} \bigl\|h(x) \bigr\|^2 F_\BX(d\mathbf {x})
= \bigl\|S(h)\bigr\|_\cH^2.
\end{eqnarray}
We shall use these equalities heavily in what follows.

Note that for every $\bs\in\bbr^d$, the fixed $\bs$ covariance
function $R_{\bs}=R(\cdot,\bs)$ is in~$\cH$, and for every $w_H\in
\cH$, and $\bt\in\bbr^d$, $w_H(\bt)=(w_H,R_\bt)_\cH$, meaning
that the
coordinate projections are continuous operations on the RKHS. This is
also the \textit{reproducing property} of the RKHS.
Note also that the quadruple $(C(T), \cH, w,
\mu_\BX)$ in general,
or $(C_0, \cH, S,\mu_\BX)$ in the stationary case, is a Wiener
quadruple in the sense of Section 3.4 in \cite{deuschelstroock1989}.

In the sequel we will use the notation $M^+(E)$ [resp., $M_1^+(E)$] for
the collection
of all Borel finite (resp., probability) measures on a topological
space~$E$.

\section{The basic large deviations result}
\label{secLD}

We start with a large deviation result for the
probability $\Psi_{\ba,\bb}$ there exists a path between $\ba$
and $\bb$
wholly within a connected component of an excursion set.

\begin{theorem} \label{tLDPpath}
\textup{(i)} Let $\BX= (X(\bt), \bt\in T)$ be a continuous Gaussian
random field on a compact set $T\subset\bbr^d$. Then
%
\begin{equation}
\label{eLDgenlimit} \lim_{u\to\infty} \frac{1}{u^2} \log
\Psi_{\ba,\bb}(u) = -\frac12 \CXab,
\end{equation}
where
%
\begin{eqnarray}
\CXab&\definedas&\inf \bigl\{ EH^2 \dvtx H\in{\mathcal L},
\mbox{ and, for some } \xi\in{\cal P}(\ba ,\bb),
\nonumber
\\[-8pt]
\\[-8pt]
\nonumber
&& \hspace*{26pt}\xi(v)\in T \mbox{ and } w_H \bigl( \xi(v) \bigr)>1, 0\leq v\leq1
\bigr\}.
\end{eqnarray}

\textup{(ii)} Let $\BX= (X(\bt), \bt\in\bbr^d)$ be a continuous
stationary Gaussian random field, with covariance function satisfying
%
\begin{equation}
\label{evanishcov} \limsup_{\| \bt\|\to\infty} R_\BX(\bt) \leq0 .
\end{equation}
Then
%
\begin{eqnarray}
\label{eLDlimit} \lim_{u\to\infty} \frac{1}{u^2} \log\Psi_\ba(u)
= -\frac12\CXa,
\end{eqnarray}
where
\begin{eqnarray*}
\CXa&\definedas& 
\inf \biggl\{ \int_{\bbr^d} \bigl\| h(
\mathbf{x})\bigr\|^2 F_\BX(d\mathbf{x})\dvtx \mbox{ for some }
\xi\in{\cal P}(\bnull,\ba)
\nonumber
\\[-8pt]
\\[-8pt]
\nonumber
&&\hspace*{28pt} \int_{\bbr^d} e^{i(\xi(v),\mathbf{x})} \bar h(\mathbf {x})
F_\BX(d\mathbf{x})>1, 0\leq v\leq1 \biggr\} .
\end{eqnarray*}
\end{theorem}
\begin{pf}
We start with putting our problem into the large deviation setup for
Gaussian measures of \cite{deuschelstroock1989}. We will use the
language of part (i) of the theorem, but the setup for part (ii) is
completely parallel. Observe that for $u>0$
\[
\Psi_{\ba,\bb}(u) = P \bigl( u^{-1}\BX\in A \bigr) ,
\]
where $A$ is the open subset of $C(T)$ given by
\[
A \equiv A_{\ba,\bb} \definedby \bigl\{ \bomega\in C(T)\dvtx \exists \xi
\in{\cal P}(\ba,\bb) \mbox{ such that } \omega \bigl(\xi(v) \bigr)>1, 0\leq v\leq1
\bigr\}.
\]
Therefore, by Theorem 3.4.5 in \cite{deuschelstroock1989}, we
conclude that
%
\begin{eqnarray}
\label{eLDestimate} -\inf_{\bomega\in A} I(\bomega) &\leq&\liminf_{u\to\infty}
\frac{1}{u^2} \log\Psi_{\ba,\bb}(u) \leq\limsup_{u\to\infty}
\frac{1}{u^2} \log\Psi_{\ba,\bb}(u)
\nonumber
\\[-8pt]
\\[-8pt]
\nonumber
& \leq&-\inf_{\bomega\in\bar A} I(\bomega)
\end{eqnarray}
for the rate function $I$ which, by Theorem 3.4.12 of
\cite{deuschelstroock1989}, can be written as
%
\begin{equation}
\label{eratefunction} I(\bomega) = \cases{ %
\tfrac12 \|
\bomega\|_\cH^2, & \quad $\mbox{if $\bomega\in\cH$},$
\vspace*{2pt}\cr
\infty,& \quad$\mbox{if $\bomega\notin\cH$,}$}
\end{equation}
for $\bomega\in C(T)$. Then \eqref{eLDestimate} already
proves the lower limit statement
\[
\liminf_{u\to\infty} \frac{1}{u^2} \log\Psi_{\ba,\bb}(u) \geq -
\frac12\CXab,
\]
valid for both parts of the theorem. Therefore, it remains to prove
the matching upper limit. Here the argument is more involved in part
(ii) of
the theorem, since noncompactness of the domain of the field requires
us to rule out the possibility of increasingly long ranging paths. We present
the argument in this case. The proof for part (i) is similar, and
easier (since we do not have to worry about paths which ``escape to
infinity'' as in
the following).

As is common with large deviation arguments, although we know that
$A=A^\circ\neq\bar A$,
this is not per se important. All that we need show is that the $\omega
$ in the
set difference $\bar A\setminus A$ do not contribute to the infimum on
the far right of
\eqref{eLDestimate}.

We start by checking that
%
\begin{equation}
\label{ebarA} \bar A \subseteq \biggl(\bigcap_{0<\delta<1}
(1-\delta)A \biggr) \cup \biggl( \bigcap_{0<\delta<1} (1-
\delta)A_1 \biggr)
\end{equation}
(in the sense of the usual multiplication of a set of functions by a
real number), where $A_1\subset C_0$ is given by
\begin{eqnarray*}
A_1 &=& \bigl\{ \bomega\in C_0\dvtx \mbox{ for every
$r>0$ there is $\bt\in\bbr^d$ with $\|\bt\|\geq r$ }
\\
&&\hspace*{6pt} \mbox{and a path $\xi\in{\cal P}(\bnull,\bt)$ such that $\bomega \bigl(\xi(v)
\bigr)>1, 0\leq v\leq1$} \bigr\}.
\end{eqnarray*}
To see this, let $\bomega\in\bar A$, so that there is a sequence
$\bomega_n\in A, n=1,2,\ldots,$ with $\bomega_n\to\bomega$ in $C_0$.
Suppose first that there is $r>0$ such
that for a subsequence $n_k\uparrow\infty$, for each $k=1,2,\ldots$
there is a path $\xi_k\in{\cal P}(\bnull, \ba)$ satisfying
$\|\xi_k(v)\|\leq r$ and $\bomega_{n_k}(\xi_k(v))>1, 0\leq v\leq
1$. Given $0<\delta<1$, choose $k$ so large that
\[
\|\bomega_{n_k}-\bomega\|_{C_0}\leq\delta/ (1+r) .
\]
Then for every $\bt\in\bbr^d$ with $\|\bt\|\leq r$ we have
$|\bomega_{n_k}(\bt)-\bomega(\bt)| \leq\delta$, so that
$\bomega(\xi_k(v))>1-\delta$ for $0\leq v\leq1$, and
$\bomega\in(1-\delta)A$.

Alternatively, suppose that such an $r>0$ does not exist. Then for every
$r>0$, for all but finitely many $n$, there is a path $\xi_n\in{\cal
P}(\bnull, \ba)$, going through a point $\bt_n$ with $\|\bt_n\|= r$,
lying within the ball of radius $r$ centered at the origin prior to
hitting the point $\bt_n$, and such that
$\bomega_n(\xi_n(v))>1, 0\leq v\leq1$. Given $r>0$ and
$0<\delta<1$, choose $n$ outside of the above exceptional finite set,
and so large that
\[
\|\bomega_{n}-\bomega\|_{C_0}\leq\delta/ (1+r) .
\]
As before, we conclude that there is a path connecting $\bnull$ and
$\bt_n$ such that the function $\bomega$ takes values above $1-\delta$
along this path. Therefore, $\bomega\in(1-\delta)A_1$, and so we have
shown \eqref{ebarA}.

Now note that since
\[
\inf_{\bomega\in(1-\delta)A} I(\bomega) =(1-\delta)^2 \inf_{\bomega\in A} I(
\bomega)
\]
for any $0<\delta<1$, the upper limit part in \eqref{eLDlimit}, and
so the result,
will follow from~\eqref{ebarA} once we check that
$ I(\bomega) =\infty$ for any $\omega\in A_1$, which we establish by showing
that $A_1\cap\cH=\varnothing$.

Suppose that, to the contrary,
there is a $\bomega= S(h)\in A_1$ for some $h\in\cH$. Fix an arbitrary
$\vep>0$. Assumption \eqref{evanishcov} guarantees the existence of
a $r_\vep>0$ such that $R_\BX(\bt)\leq\vep$ if $\|\bt\|\geq
r_\vep$. By the
definition of $A_1$, for every $n=1,2,\ldots$ there is $\bt_n$
with $\|\bt_n\|=nr_\vep$ and a path $\xi$ connecting $\bnull$ and
$\bt_n$ such that $\bomega(\xi(v))>1, 0\leq v\leq1$. We can
choose $0<v_1<\cdots<v_n\leq1$ such that $\|\xi(v_j)\|=jr_\vep$ for
$j=1,\ldots, n$. Then
\begin{eqnarray*}
1 &<& \frac1n \sum_{j=1}^n \bomega \bigl(
\xi(v_j) \bigr) = \int_{\bbr^d} \Biggl( \frac1n
\sum_{j=1}^n e^{i(\xi(v_j),\mathbf{x})} \Biggr) \bar
h( \mathbf{x}) F_\BX (d\mathbf{x})
\\
&\leq& \Biggl\llVert \frac1n \sum_{j=1}^n
e^{i(\xi(v_j),\cdot)} \Biggr\rrVert_{L^2(F_\BX)} \| h\|_{L^2(F_\BX)} .
\end{eqnarray*}
However,
\begin{eqnarray*}
\Biggl\llVert \frac1n \sum_{j=1}^n
e^{i(\xi(v_j),\cdot)} \Biggr\rrVert_{L^2(F_\BX)}^2 &=& \frac{1}{n^2}
\Biggl( nR_\BX(0) + 2\sum_{j_1=1}^{n-1}
\sum_{j_2=j_1+1}^n R_\BX \bigl(
\xi(v_{j_1})-\xi(v_{j_2}) \bigr) \Biggr)
\\
&\leq& \frac1n R_\BX(0) + \vep,
\end{eqnarray*}
so that
\[
\| h\|^2_{L^2(F_\BX)}> \frac{1}{\frac1n R_\BX(0) + \vep} .
\]
Sending first $n\to\infty$ and then $\vep\to0$, we obtain $\|
h\|_{L^2(F_\BX)}=\infty$, which is impossible.

This contradiction proves the rightmost inequality in
\eqref{eLDlimit} and so we are done.
\end{pf}

Theorem \ref{tLDPpath} describes the logarithmic asymptotic
of the path existence probability $\Psi_{\ba,\bb}$ in terms of a
solution to an optimization problem in the Hilbert space. The next
result contains the dual version of this optimization problem and
relates $\Psi_{\ba,\bb}$ to the problem of finding a path
of minimal capacity between $\ba$ and $\bb$.

\begin{theorem} \label{tlargeenergypath}
\textup{(i)} Let $\BX= (X(\bt), \bt\in T)$ be a continuous Gaussian
random field on a compact set $T\subset\bbr^d$. Then
%
\begin{eqnarray}
\label{ecapGenlimit} && \lim_{u\to\infty} \frac{1}{u^2} \log
\Psi_{\ba,\bb}(u)
\nonumber\\
&&\qquad = -\frac12\CXab
\\
&&\qquad = -\frac12 \biggl[ \sup_{\xi\in\cP(\ba,\bb)} \min_{\mu\in
M_1^+([0,1])} \int
_0^1\int_0^1
R_\BX \bigl( \xi(u),\xi(v) \bigr) \mu(du) \mu(dv)
\biggr]^{-1}.
\nonumber
\end{eqnarray}
%

\textup{(ii)} Let $\BX= (X(\bt), \bt\in\bbr^d)$ be a continuous
stationary Gaussian random field, with covariance function satisfying
\eqref{evanishcov}. Then
%
\begin{eqnarray}
\label{ecaplimit} &&\lim_{u\to\infty} \frac{1}{u^2} \log
\Psi_\ba(u)
\nonumber\\
&&\qquad = -\frac12\CXa
\\
&&\qquad = -\frac12 \biggl[ \sup_{\xi\in\cP(\bnull,\ba)} \min_{\mu\in
M_1^+([0,1])} \int
_0^1\int_0^1
R_\BX \bigl( \xi(u)-\xi(v) \bigr) \mu(du) \mu(dv)
\biggr]^{-1}.
\nonumber
\end{eqnarray}
%
\end{theorem}
Note that the space $M_1^+([0,1])$
is weakly compact, and the covariance function $R_\BX$ is
continuous. Therefore, for a fixed path $\xi$, the function
\[
\mu\to\int_0^1\int_0^1
R_\BX \bigl( \xi(u),\xi(v) \bigr) \mu(du) \mu(dv)
\]
is weakly continuous on compacts. Hence, it achieves its infimum,
and it is legitimate to write ``min'' in \eqref{ecapGenlimit} and in
\eqref{ecaplimit}.
\begin{pf*}{Proof of Theorem \ref{tlargeenergypath}}
The proofs of the two parts are only notationally different, so we
will suffice with a proof for part (i) only.
We use the Lagrange duality approach of Section 8.6 in
\cite{luenberger1969}. Writing
\[
\CXab= \inf_{\xi\in\cP(\ba,\bb)} \CXbxi,
\]
where, for $\xi\in\cP(\ba,\bb)$,
%
\begin{eqnarray}
\label{eprimalfeasible} \CXbxi\definedas\inf \bigl\{ EH^2 \dvtx H
\in{\mathcal L} \mbox{ and } w_H \bigl( \xi(v) \bigr)>1, 0\leq v\leq1
\bigr\},
\end{eqnarray}
we see that it is enough to prove that for every $\xi\in
\cP(\ba,\bb)$,
%
\begin{eqnarray}
\label{epathdual}&& \CXbxi
\nonumber
\\[-8pt]
\\[-8pt]
\nonumber
&&\qquad= \biggl[ \min_{\mu\in
M_1^+([0,1])} \int_0^1
\int_0^1 R_\BX \bigl( \xi(u),
\xi(v) \bigr) \mu(du) \mu(dv) \biggr]^{-1}.
\end{eqnarray}
To this end, let ${\tt Z}=C([0,1]$. Then ${\tt P} \definedas\{ z\in
{\tt
Z}\dvtx z(v)\geq0, 0\leq v\leq1\}$ is a closed convex cone in
${\tt Z}$. Its dual cone ${\tt P}^\oplus\subset{\tt Z}^*$ [defined
as the collection of $z^*\in{\tt Z}^*$ such that $z^*(z)\geq0$ for
all $z\in{\tt P} $] can be
naturally identified with $M^+([0,1])$. Fix $\xi\in\cP(\ba,\bb)$, and
define $G\dvtx {\mathcal L}
\to{\tt Z}$ by
\[
G(H) = G_\xi(H)\definedby \bigl( 1-w_H \bigl( \xi(v)
\bigr), 0\leq v\leq1 \bigr).
\]
Then $G$ is, clearly, a convex mapping. We can also write
%
\begin{equation}
\label{eprimal} \bigl( \CXbxi \bigr)^{1/2} =\inf \bigl\{ \bigl(EH^2
\bigr)^{1/2}\dvtx H\in {\mathcal L}, G(H)\in -{\tt P} \bigr\},
\end{equation}
and so our task now is to show that \eqref{eprimal} implies \eqref
{epathdual}.

Suppose first that the feasible set in the optimization problem
\eqref{eprimalfeasible} is not empty. Then there is $H\in
{\mathcal L}$ such that $G(H)$ belongs to the interior of the cone
$-{\tt P}$, so by Theorem 1, page 224 of \cite{luenberger1969}, we
conclude that
%
\begin{eqnarray}
\label{edualconst}&& \bigl( \CXbxi \bigr)^{1/2}
\nonumber
\\[-8pt]
\\[-8pt]
\nonumber
&&\qquad=\max_{\mu\in M^+([0,1])}
\inf_{H\in{\mathcal L}} \biggl[ \bigl( EH^2 \bigr)^{1/2}+ \int
_0^1 G(H) (v) \mu(dv) \biggr] ,
\end{eqnarray}
and we may use ``max'' instead of ``sup'' because an optimal $\mu\in
M^+([0,1])$ exists. For a fixed
$\mu\in M^+([0,1])$ with total mass $\|\mu\|$, we let $\hat\mu=\mu
/\|\mu\|\in M^+_1([0,1])$. Then
%
\begin{eqnarray}
\label{elagrange}
&& \inf_{H\in{\mathcal L}} \biggl[ \bigl( EH^2
\bigr)^{1/2} + \int_0^1 G(H) (v)
\mu(dv) \biggr]
\nonumber\\
&&\qquad = \|\mu\| + \inf_{H\in{\mathcal L}} \biggl[ \bigl( EH^2
\bigr)^{1/2}- \|\mu\|\int_0^1
w_H \bigl( \xi(v) \bigr) \hat\mu(dv) \biggr]
\\
&& \qquad= \|\mu\| + \inf_{a\geq0} a \biggl[ 1- \|\mu\|
\sup_{H\in{\mathcal
L}\dvtx EH^2=1} \int_0^1 w_H
\bigl( \xi(v) \bigr) \hat\mu(dv) \biggr]
\nonumber\\
&&\qquad = \cases{ %
 -\infty,& \quad $\mbox{if } \displaystyle\|\mu\|>
\biggl[\sup_{H\in{\mathcal L}\dvtx EH^2=1} \int_0^1
w_H \bigl( \xi(v) \bigr) \hat\mu(dv) \biggr]^{-1},$
\vspace*{2pt}\cr
\|\mu\|, &\quad  $\mbox{if } \displaystyle\|\mu\|\leq \biggl[ \sup_{H\in{\mathcal
L}\dvtx EH^2=1} \int
_0^1 w_H \bigl( \xi(v) \bigr) \hat
\mu(dv) \biggr]^{-1}.$}
\nonumber
\end{eqnarray}
Therefore,
\[
\bigl( \CXbxi \bigr)^{1/2} = \biggl[ \inf_{\mu\in M_1^+([0,1])} \sup_{H\in{\mathcal L}\dvtx EH^2=1} \int
_0^1 w_H \bigl( \xi(v) \bigr)
\mu(dv) \biggr]^{-1} ,
\]
and \eqref{epathdual} follows, since by the reproducing property of
the RKHS, for every
$\mu\in M_1^+([0,1])$,
\begin{eqnarray*}
\sup_{H\in{\mathcal L}\dvtx EH^2=1} \int_0^1 w_H
\bigl( \xi(v) \bigr) \mu(dv) &=& \sup_{H\in{\mathcal
L}\dvtx EH^2=1} \int_0^1
\bigl( w_H, R_X \bigl(\xi(v),\cdot \bigr)
\bigr)_\cH\mu(dv)
\\
&=& \sup_{w\in\cH\dvtx  \| w\|_\cH=1} \biggl( w, \int_0^1
R_X \bigl(\xi (v),\cdot \bigr) \mu(dv) \biggr)_\cH
\\
&=& \biggl( \int_0^1\int_0^1
R_\BX \bigl( \xi(u),\xi(v) \bigr) \mu(du) \mu(dv)
\biggr)^{1/2}.
\end{eqnarray*}
In the last step we have used the fact that
\[
w_\mu\definedas\int_0^1
R_X \bigl(\xi(v),\cdot \bigr) \mu(dv) \in\cH,
\]
so the supremum of the inner product is achieved at
$w = w_\mu/\| w_\mu\|_{\cal H}$, and
\[
\| w_\mu\| = \biggl( \int_0^1\int
_0^1 R_\BX \bigl( \xi(u),\xi (v)
\bigr) \mu(du) \mu(dv) \biggr)^{1/2} .
\]

This establishes \eqref{epathdual} for the case that the feasible
set in
\eqref{eprimalfeasible} is not empty. We now turn to the case in
which this set is, indeed,
empty. This will complete the proof of the theorem.
In this case \eqref{epathdual}
reduces to the statement
%
\begin{equation}
\label{eIstar} I_* \definedas\min_{\mu\in M_1^+([0,1])} \int_0^1
\int_0^1 R_\BX \bigl( \xi(u),
\xi(v) \bigr) \mu(du) \mu(dv) =0 .
\end{equation}
Suppose that, to the contrary, $I_*>0$. Let $\mu_0\in M_1^+([0,1])$
achieve the minimum value in the integral defining $I_*$. Consider the
continuous real-valued function
\[
W(u) =\int_0^1 R_\BX \bigl(
\xi(u),\xi(v) \bigr) \mu_0(dv),\qquad  0\leq u\leq1.
\]
If this function never vanishes, then, by continuity and compactness,
it is bounded away from zero, so a sufficiently large in absolute
value multiple of the random variable in $\mathcal L$ given by
\[
H = \int_0^1 X \bigl( \xi(v) \bigr)
\mu_0(dv)
\]
is feasible for the optimization problem
\eqref{eprimalfeasible}, contradicting the assumption that the set
of feasible solutions is empty.

Hence, there is $u_0\in[0,1]$ such
that $W(u_0) =0$. For $0<\vep<1$ define a probability measure in
$M_1^+([0,1])$ by
\[
\mu_\vep= (1-\vep)\mu_0 + \vep\delta_{u_0} ,
\]
where $\delta_a$ denotes the point mass at $a$. Note that
\begin{eqnarray*}
I(\vep)&\definedas& \int_0^1\int
_0^1 R_\BX \bigl( \xi(u),\xi (v)
\bigr) \mu_\vep(du) \mu_\vep(dv)
\\
&=& (1-\vep)^2 I_* + 2\vep(1-\vep) W(u_0) +
\vep^2R_\BX(0)
\\
&=& (1-\vep)^2 I_*+ \vep^2R_\BX(0) .
\end{eqnarray*}
Since $I_*$ was assumed to be positive, we see that
\[
\frac{d}{d\vep}I(\vep)\Big|_{\vep=0}<0 ,
\]
which contradicts the minimality of $I_*$. This proves
\eqref{eIstar} and so the theorem.
\end{pf*}

Observe that an alternative way of stating the result of Theorem
\ref{tlargeenergypath} is
%
\begin{equation}
\label{eenergyalt} \CXab= \inf_{\xi\in\cP(\ba,\bb)} \biggl[ \min_{\mu\in
M_1^+(\xi)} \int
_\xi\int_\xi R_\BX(\bt,
\bs) \mu(d\bt) \mu (d\bs) \biggr]^{-1} ,
\end{equation}
where $M_1^+(\xi)$ is the set of all probability measures in $\bbr^d$
supported by the path $\xi$ (strictly speaking, by the compact image
of the interval $[0,1]$ under $\xi$). For a fixed path $\xi\in
\cP(\ba,\bb)$, the quantity
%
\begin{equation}
\label{emeasurepath} \CXbxi= \biggl[ \min_{\mu\in
M_1^+(\xi)} \int
_\xi\int_\xi R_\BX(\bt,
\bs) \mu(d\bt) \mu (d\bs) \biggr]^{-1}
\end{equation}
is known as \textit{the capacity of the path $\xi$ with respect to
the kernel} $R_\BX$;
see \cite{fuglede1960}. Therefore, we can treat the problem of solving
\eqref{eenergyalt} as one of finding a path between the points $\ba
$ and
$\bb$ of minimal capacity.

\section{Fixed paths and measures of minimal
energy} \label{secsmallenergy}

The dual formulation \eqref{eLDlimit} of the optimization problem
required to find the asymptotics of the path existence
probability $\Psi_{\ba,\bb}(u)$ involves solving fixed path $\xi$
optimization problems~\eqref{eprimalfeasible} or~\eqref{epathdual}. For a fixed path we have the following version of
Theorems \ref{tLDPpath} and \ref{tlargeenergypath}.
%
\begin{theorem} \label{tfixedpath}
\textup{(i)} For a $\xi\in\cP(\ba,\bb)$ let
\[
\Psi_{\ba,\bb}(u; \xi) = P \bigl( X \bigl(\xi(v) \bigr)>u, 0\leq v\leq1
\bigr) .
\]
Then
%
\begin{equation}
\label{epsiLD} \lim_{u\to\infty} \frac{1}{u^2} \log\Psi_{\ba,\bb}(u;
\xi) = - \frac12\CXbxi.
\end{equation}
\begin{longlist}[(iii)]
\item[(ii)] The primal problem \eqref{eprimalfeasible} can be
rewritten in the form
%
\begin{eqnarray}
\label{eprimalpath} \CXbxi= \inf \bigl\{ EH^2\dvtx H\in{\mathcal L},
E \bigl[ X \bigl( \xi(v) \bigr)H \bigr]\geq1, 0\leq v\leq1 \bigr\}.
\end{eqnarray}
Further, if the feasible set in \eqref{eprimalpath} is nonempty,
then the infimum in \eqref{eprimalpath} is achieved at a unique
$H_\xi\in
{\mathcal L}$.

\item[(iii)] The set $\Wxi$ of $\mu\in M_1^+([0,1])$ over which the
minimum in the dual problem~\eqref{epathdual} is achieved is a
weakly compact
convex subset\vadjust{\goodbreak} of $M_1^+([0,1])$. Furthermore, if the feasible set in
\eqref{eprimalpath} is nonempty, then,
for every $\mu\in\Wxi$,
%
\begin{equation}
\label{ecomplslack} \mu \bigl( \bigl\{ 0\leq v\leq1\dvtx E \bigl[ X \bigl(
\xi(v) \bigr)H_\xi \bigr] >1 \bigr\} \bigr) = 0 .
\end{equation}

\item[(iv)] Suppose that the feasible set in \eqref{eprimalpath} is
nonempty. Then for every $\vep>0$,
%
\begin{equation}
\label{emostlikely} P \biggl( \sup_{0\leq v\leq1} \biggl\llvert \frac1u X \bigl(
\xi(v) \bigr) - x_\xi(v) \biggr\rrvert \geq\vep \Big| X \bigl(\xi(v)
\bigr)>u, 0\leq v\leq 1 \biggr)\to0
\end{equation}
as $u\to\infty$. Here
%
\begin{equation}
\label{elikelypath} x_\xi(v) =E \bigl[ X \bigl( \xi(v)
\bigr)H_\xi \bigr],\qquad 0\leq v\leq1 .
\end{equation}
\end{longlist}
\end{theorem}
The probability measures $\mu\in\Wxi$ are called \textit{capacitary
measures}, or \textit{measures of minimal energy}; see \cite{fuglede1960}.
\begin{pf*}{Proof of Theorem \ref{tfixedpath}}
Part (i) of the theorem can be proved in the same way as Theorem
\ref{tLDPpath}. The fact that the primal formulations
\eqref{eprimalfeasible} and \eqref{eprimalpath} are equivalent is
an immediate consequence of the definition of $w_H$. Suppose now that
the feasible set in \eqref{eprimalpath} is
nonempty, and let $H_n\in{\mathcal L}$, $n=1,2,\ldots$ be a sequence of
feasible solutions such that $EH_n^2\to\CXbxi$. The weak
compactness of the unit ball in $\mathcal L$ shows that this sequence
has a subsequential weak limit $H_\xi$ with $EH_\xi^2 =
\CXbxi$. Since the set of feasible solutions is weakly closed, $H_\xi$
is feasible. The uniqueness of the optimal solution to
\eqref{eprimalpath} follows from convexity of the norm.

Convexity and weak compactness of the set $\Wxi$ follow from the
nonnegative definiteness and continuity of $R_\BX$; see, for example,
Remark 2,
page 160, in \cite{fuglede1960}. The statement \eqref{ecomplslack} is
a part of the relation between the dual and primal optimal solutions;
see Theorem 1, page 224, in \cite{luenberger1969}.

For part (iv) of the theorem, note that by the Gaussian large deviation
principle of Theorem 3.4.5 in \cite{deuschelstroock1989},
%
\begin{eqnarray}\label{emoreconst}
\nonumber\qquad
&&\limsup_{u\to\infty} \frac{1}{u^2} \log P \biggl( X \bigl(
\xi(v) \bigr)>u, 0\leq v\leq 1, \sup_{0\leq v\leq1} \biggl\llvert \frac1u X \bigl(
\xi(v) \bigr) - x_\xi(v) \biggr\rrvert \geq\vep \biggr)
\\
 &&\qquad \leq-\frac12 \inf \Bigl\{ EH^2, H\in{\mathcal
L}\dvtx E \bigl[ X \bigl( \xi(v) \bigr)H \bigr]\geq1, 0\leq v\leq1,
\\
&&\hspace*{25pt}\qquad\qquad \sup_{0\leq v\leq1} \bigl\llvert E \bigl[ X \bigl( \xi(v) \bigr)H \bigr] - E
\bigl[ X \bigl( \xi(v) \bigr)H_\xi \bigr] \bigr\rrvert \geq\vep \Bigr
\} .
\nonumber
\end{eqnarray}
Therefore, the statement \eqref{emostlikely} will follow from Parts
(i) and (ii) of the theorem once we prove that the infimum in
\eqref{emoreconst} is strictly larger than $\CXbxi$. Suppose that,
to the contrary, the two infima are equal. By the weak
compactness of the unit ball in $\mathcal L$ and the fact that the
feasible set in \eqref{emoreconst} is weakly closed, this would
imply existence of $H_*$ feasible for \eqref{emoreconst} such that
$EH_*^2= EH_\xi^2$. Since $H_\xi$ is not
feasible for \eqref{emoreconst},\vadjust{\goodbreak} we know that $H_*\not=
H_\xi$. Since $H_*$ is feasible for \eqref{eprimalpath}, we have
obtained a contradiction to the uniqueness of $H_\xi$ proved
above. This completes the proof of the theorem.~%
\end{pf*}

\begin{remark} \label{rkinterpret}
Theorem \ref{tfixedpath} has the following important
interpretation. Assuming that the feasible set in \eqref
{eprimalpath} is
nonempty, part (iv) of the theorem implies that the nonrandom function
$x_\xi$ in \eqref{elikelypath} is
the most likely choice for the normalized sample path
$u^{-1}X(\xi(v)), 0\leq v\leq1$ along $\xi$, given that
$\{X(\xi(v))>u, 0\leq v\leq
1\}$. Part (iii) of the theorem implies that the values of the
random field along the path $\xi$ have to
(nearly) touch the level $u$ at the points of the support of
any measure of minimal energy. In other words, the sample
path needs to be ``supported,'' or ``held,'' at the level $u$ at
the points of the support in order to achieve the highest
probability of exceeding the high level $u$ along the entire path
$\xi$. We will see explicit examples of how this works in the
following section,
when we more closely investigate the one-dimensional case.
\end{remark}

The duality relation of the optimization problems
\eqref{eprimalpath} and \eqref{epathdual} immediately provides
upper and lower bounds on $\CXbxi$ of the form
%
\begin{equation}
\label{eboundsexp} \biggl[ \int_0^1\int
_0^1 R_\BX \bigl( \xi(u),\xi(v)
\bigr) \mu(du) \mu(dv) \biggr]^{-1} \leq\CXbxi\leq EH^2
\end{equation}
for any $\mu\in M_1^+([0,1])$ and any $H\in{\mathcal L}$ feasible for
\eqref{eprimalpath}. In particular, if
%
\begin{equation}
\label{ecriterionopt} \biggl[ \int_0^1\int
_0^1 R_\BX \bigl( \xi(u) \xi(v)
\bigr) \mu(du) \mu(dv) \biggr]^{-1} =EH^2
\end{equation}
for some $\mu$ and $H$ as above, then $\mu\in\Wxi$, $H=H_\xi$, and
the common value in \eqref{ecriterionopt} is equal to $\CXbxi$.

Finding a measure of minimal energy, $\mu\in
\Wxi$, is, in general, a difficult problem. The following theorem
includes a characterization of these measures.
%
\begin{theorem} \label{toptimalmchararct}
Assume that the feasible set in \eqref{eprimalpath} is
nonempty.
\begin{longlist}[(ii)]
\item[(i)] For every $\mu\in\Wxi$ we have
\[
H_\xi= \CXbxi\int_0^1 X \bigl(
\xi(v) \bigr) \mu(dv)
\]
with probability 1.

\item[(ii)] A probability measure $\mu\in M^+_1([0,1])$ is a measure
of minimal energy (i.e., $\mu\in\Wxi$) if and only if
%
\begin{eqnarray}
\label{ehighestenergy} &&\min_{0\leq v\leq1} \int_0^1
R_\BX \bigl( \xi(u),\xi(v) \bigr) \mu(du)
\nonumber
\\[-8pt]
\\[-8pt]
\nonumber
&& \qquad= \int_0^1\int_0^1
R_\BX \bigl( \xi(u_1),\xi(u_2) \bigr)
\mu(du_1) \mu(du_2) >0 .
\end{eqnarray}
\end{longlist}
\end{theorem}
Note that part (ii) of the theorem also says that the integral in the
left-hand side of \eqref{ehighestenergy} is equal to the double
integral in its right-hand side for $\mu$-almost every $0\leq v\leq
1$.
\begin{pf*}{Proof of Theorem \ref{toptimalmchararct}}
For part (i), let $\mu\in\Wxi$. The calculations following
the maximization problem \eqref{edualconst} show that $\mu_\xi=
\CXbxi^{1/2}\mu$ is an optimal measure for that problem. It follows
from Theorem 1, page 224, in \cite{luenberger1969} that $H_\xi$ solves
the minimization problem in \eqref{elagrange}, when any measure in
$M^+([0,1])$ optimal for
\eqref{edualconst} is used. Using the
measure $\mu_\xi$, we see that
\[
H_\xi=a \int_0^1 X \bigl( \xi(v)
\bigr) \mu_\xi(dv) = a \CXbxi^{1/2} \int_0^1
X \bigl( \xi(v) \bigr) \mu(dv)
\]
for some $a>0$. Testing all random variables of the type
\[
H_\xi=b \int_0^1 X \bigl( \xi(v)
\bigr) \mu(dv),\qquad b>0,
\]
in \eqref{eprimalpath} leads to the conclusion that
\[
b = \biggl[ \min_{0\leq v\leq1} \int_0^1
R_\BX \bigl( \xi(u),\xi(v) \bigr) \mu(du) \biggr]^{-1}.
\]
The fact that $b=\CXbxi$ follows now from the optimality of $\mu$ and
the general properties of measures of minimal energy for bounded symmetric
kernels; see, for example, Theorem 2.4 in \cite{fuglede1960}.

We now prove part (ii). Suppose first that $\mu$ satisfies
\eqref{ehighestenergy}, and define
$H\in{\mathcal L}$ by
\[
H= \bigl(K(\mu) \bigr)^{-1} \int_0^1
X \bigl( \xi(v) \bigr) \mu(dv),
\]
where $K(\mu)$ is the double integral in
the right-hand side of \eqref{ehighestenergy}. Note that for any
$0\leq v\leq1$,
\[
E \bigl[ X \bigl( \xi(v) \bigr)H \bigr] = \bigl(K(\mu) \bigr)^{-1} \int
_0^1 R_\BX \bigl( \xi(u),\xi(v)
\bigr) \mu(du) \geq1
\]
by \eqref{ehighestenergy}. Therefore, $H$ is feasible for
\eqref{eprimalpath}. However,
\[
EH^2 = \frac{1}{K(\mu)^2} E \biggl( \int_0^1
X \bigl( \xi(v) \bigr) \mu(dv) \biggr)^2 = \frac{1}{K(\mu)} ,
\]
so that $\mu$ and $H$ satisfy the relation
\eqref{ecriterionopt}. Hence, $\mu\in\Wxi$ (and $H=H_\xi$).

In the opposite direction, if $\mu\in\Wxi$, then the equality
in \eqref{ehighestenergy} is a general property of measures of
minimal energy for bounded symmetric kernels, as in the proof of part
(i). The fact that the equal terms in \eqref{ehighestenergy} are
positive follows from the fact that the feasible set in \eqref
{eprimalpath} is
nonempty, so $\CXbxi<\infty$.
\end{pf*}

\begin{remark} \label{rkstationery}
If the random field $\BX$ is stationary, then the results of this
section can be restated in the language used in Section \ref{secLD}
in the stationary case. In particular, the primal problem
\eqref{eprimalpath} becomes
\begin{eqnarray*}
\CXxi = \inf \biggl\{ \int_{\bbr^d} \bigl\| h(\mathbf{x})
\bigr\|^2 F_\BX(d\mathbf{x})\dvtx %
%
\int_{\bbr^d} e^{i(\xi(v),\mathbf{x})}
\bar h(\mathbf{x}) F_\BX (d\mathbf{x})>1, 0\leq v\leq1 \biggr\},
\end{eqnarray*}
while the optimal solution of the primal problem in part (i) of
Theorem \ref{toptimalmchararct} becomes
\[
h_\xi(\mathbf{x})= \CXxi\int_0^1
e^{i(\xi(v),\mathbf{x})} \mu (dv)\qquad \mbox{for $F_\BX$-almost all } \mathbf{x}\in
\bbr^d
\]
and for any $\mu\in\Wxi$. This relation can also be restated in terms on
the measures
supported by the (image of) path $\xi$ instead of the unit interval,
as in \eqref{emeasurepath}. If $\mu$ is an optimal measure in
\eqref{emeasurepath}, then we have
%
\begin{equation}
\label{euptonull} h_\xi(\mathbf{x})= \CXxi\int_\xi
e^{i(\bt,\mathbf{x})} \mu (d\bt) \qquad\mbox{for $F_\BX$-almost all } \mathbf{x}
\in\bbr^d .
\end{equation}
Note that the function in the right-hand side of \eqref{euptonull}
is, up to a constant, the characteristic function of
the measure $\mu$. If the support of the spectral measure $F_\BX$
happens to be the entire space $\bbr^d$, then the characteristic
functions of all optimal measures in \eqref{emeasurepath} are equal
and, hence, the uniqueness of a characteristic function shows
that, in this case [and as long as the feasible set in \eqref
{eprimalpath} is
nonempty], there is exactly one probability measure $\mu\in
M_1^+(\xi)$ of minimal energy.
\end{remark}

\begin{remark} \label{rkconstant}
An immediate conclusion of part (ii) of Theorem
\ref{toptimalmchararct} and the assumed continuity of the
covariance function
is that the function
\[
v\mapsto\int_0^1 R_\BX \bigl(
\xi(u),\xi(v) \bigr) \mu(du), \qquad 0\leq v\leq1,
\]
is constant on the support of any measure $\mu\in\Wxi$. This seems to
indicate that the support of any measure of minimal energy may
not be ``large.'' In the examples below, however, this intuition holds
only in some cases.
\end{remark}

\section{The one-dimensional case}
\label{seconedim}

In this and the following two sections we specialize to the
one-dimensional case $d=1$. Let
$a<b$. As before, we are interested in the probability
\[
\Psi_{a,b}(u)= P \bigl( X(t)>u, a\leq t\leq b \bigr) .
\]
There is, essentially, a single path between $a$ and $b$,
and the results of the previous two sections immediately specialize to
yield the following
special case. [Note that
condition \eqref{evanishcov} is superfluous in the one-dimensional
nonstationary case.]

\begin{theorem} \label{t1dimLDP}
Let $\BX$ be a continuous Gaussian process on an interval including $[a,b]$.
Then the limit
\[
-\frac12 C_\BX(a,b)\definedas\lim_{u\to\infty} \frac{1}{u^2}
\log\Psi_{a,b}(u)
\]
exists, and
%
\begin{eqnarray}
\quad && C_\BX(a,b)
\nonumber
\\
&&\qquad = \inf \bigl\{ EH^2\dvtx H\in{\mathcal L}, E \bigl[ X \bigl(
a+(b-a)v \bigr)H \bigr]\geq1, 0\leq v\leq1 \bigr\} \label{e1dimrepr}
\\
\label{e1dimreprdual} && \qquad = \biggl[ \min_{\mu\in M_1^+([0,1])} \int_0^1
\int_0^1 R_\BX \bigl(
a+(b-a)u,a+(b-a)v \bigr)  \mu(du) \mu(dv) \biggr]^{-1}.
\end{eqnarray}
If the process $\BX$ is stationary, an alternative expression for
$C_\BX(a)\definedas C_\BX(0,a)$, $a>0$, is given by
%
\begin{eqnarray}
\label{e1dimreprst}
&& C_\BX(a) = \inf \biggl\{ \int_{-\infty}^\infty
\bigl\| h(x)\bigr\|^2 F_\BX (dx) \dvtx \int_{-\infty}^\infty
e^{ivax}\bar h(x) F_\BX(dx)>1,
\nonumber
\\[-8pt]
\\[-8pt]
\nonumber
 &&\hspace*{235pt}0\leq v\leq1 \biggr\}.
\end{eqnarray}
The set $\Wab$ of $\mu\in M_1^+([0,1])$
over which the minimum in \eqref{e1dimreprdual} is achieved is a
weakly compact convex subset of $M_1^+([0,1])$. The measures in $\Wab$
are characterized by the relation
%
\begin{eqnarray}
\label{echeck1dim} &&\min_{0\leq v\leq1} \int_0^1
R_\BX \bigl( a+(b-a)u,a+(b-a)v \bigr) \mu(du)
\nonumber
\\[-8pt]
\\[-8pt]
\nonumber
&&\qquad = \int_0^1\int
_0^1 R_\BX \bigl(
a+(b-a)u_1,a+(b-a)u_2 \bigr) \mu(du_1)
\mu(du_2) .
\end{eqnarray}

Suppose, further, that the problem \eqref{e1dimrepr} has a feasible
solution. In this case the double integral in \eqref{echeck1dim} is
positive for any $\mu\in\Wab$, and the problem \eqref{e1dimrepr}
has a unique optimal solution, $H_{a,b}$. For each $\mu\in\Wab$,
\[
H_{a,b}= C_\BX(a,b) \int_0^1
X \bigl( a+(b-a)v \bigr) \mu(dv)
\]
with probability 1. In the stationary case, the problem \eqref{e1dimreprst}
has a unique optimal solution, $h_{a}$. For each $\mu\in
\Wa\definedas{\mathcal W}_{0,a}$
\[
h_a(x) =C_\BX(a) \int_0^1
e^{iavx} \mu(dv)\qquad \mbox{for $F_\BX$-almost all $-\infty<x<
\infty$.}
\]

The conditional law on $C[a,b]$ of the scaled process $u^{-1}\BX$
restricted to the interval $[a,b]$, given that $X(t)>u, a\leq t\leq
b$, converges as $u\to\infty$ to the Dirac measure at
%
\begin{equation}
\label{emostlikelya} x_{a,b}(t) = C_\BX(a,b) \int
_0^1 R_\BX \bigl( t,a+(b-a)v \bigr)
\mu(dv),\qquad  a\leq t\leq b
\end{equation}
and
\[
\mu \bigl( \bigl\{ 0\leq v\leq1\dvtx x_{a,b} \bigl( a+(b-a)v \bigr)>1
\bigr\} \bigr) = 0 .
\]

Finally, if the process $\BX$ is stationary, and the support of
$F_\BX$ is the entire real line, then
the set $\Wa$ consists of a single probability measure, $\mu_a$.
\end{theorem}

\begin{remark} \label{rksymmetry}
Suppose that the process $\BX$ is stationary. For $\mu\in
M_1^+([0,1])$ define $\hat\mu= \mu\circ
T^{-1}$ with $T\dvtx [0,1]\to[0,1]$ being the reflection map $Tx= 1-x$,
$0\leq x\leq1$. If $\mu\in\Wa$, then $\hat\mu$ satisfies conditions
\eqref{echeck1dim} because $\mu$ does, hence $\hat\mu\in\Wa$ as
well. By convexity of $\Wa$, so does the symmetric (around $x=1/2$)
probability measure $1/2(\mu+\hat\mu)$. Therefore, $\Wa$ always
contains a symmetric measure. In particular, if $\Wa$ is a singleton,
then the unique measure of minimal energy is symmetric.
\end{remark}

In the remainder of this section we concentrate on the stationary
case. We will investigate how the
probability measure $\mu_a$, the function $h_a$ and the limiting shape
$x_a \definedas x_{0,a}$ change as functions
of $a$. This will help us understand the order of magnitude of the
probability $\Psi_a(u)$ for varying lengths $a$ of the interval and,
according to part (iv) of Theorem \ref{tfixedpath}, it will tell us
the most likely shape the process $\BX$ takes when it exceeds a high
level $u$ along the entire interval $[0,a]$.

Our first result describes the situation occurring for some, but not all,
stationary Gaussian processes on short intervals.
%
\begin{proposition} \label{prcase1}
Let $\BX$ be a stationary continuous Gaussian process. Suppose that
for some $a>0$ the following condition holds:
%
\begin{equation}
\label{econd1} R_\BX(t)+R_\BX(a-t)\geq
R_\BX(0)+R_\BX(a)>0 \qquad\mbox{for all $0\leq t\leq a$.}
\end{equation}
Then a measure in $\Wa$ is given by
%
\begin{equation}
\label{efirstmu} \mu^{(1)} \definedby\tfrac12 \delta_0+
\tfrac12 \delta_1.
\end{equation}
Furthermore,
%
\begin{eqnarray}
\label{efirstC} C_\BX(a)&=& \frac{2}{R_\BX(0)+R_\BX(a)},
\\
\label{efirsth} h_a(x) &= &\frac{1+e^{iax}}{R_\BX(0)+R_\BX(a)} \qquad\mbox{for
$F_\BX$-almost all $-\infty<x<\infty$}
\end{eqnarray}
and
%
\begin{equation}
\label{efirstx} x_a(t) = \frac{R_\BX(t) + R_\BX(a-t)}{R_\BX(0)+R_\BX(a)},\qquad 0\leq t\leq a .
\end{equation}
\end{proposition}
\begin{pf}
Once we show that $\mu^{(1)}\in\Wa$, the rest of the statements will
follow from Theorem \ref{t1dimLDP}. In order to prove \eqref
{efirstmu}, we
need to check conditions \eqref{echeck1dim}. These follow
immediately from \eqref{econd1} and the fact that
\[
\int_0^1\int_0^1
R_\BX \bigl( a(u_1-u_2) \bigr)
\mu^{(1)}(du_1) \mu^{(1)}(du_2) =
\frac12 R_\BX(0)+\frac12 R_\BX(a) ,
\]
while for $0\leq v\leq1$,
\[
\int_0^1 R_\BX \bigl( a(u-v)
\bigr) \mu^{(1)} (du) = \frac12 R_\BX(av)+\frac12
R_\BX(a-av) .
\]
\upqed\end{pf}
%
\begin{remark} \label{rkcase1}
Note that a sufficient (but not necessary) condition for
\eqref{econd1} is concavity of the covariance function $R_\BX$ on
the interval $[0,a]$. Indeed, for a concave covariance function the derivative
exists apart from a countable set of points and is
monotone. Therefore,
\[
R_\BX(t) - R_\BX(0) = \int_0^t
R_\BX^\prime(s) \,ds \geq \int_0^t
R_\BX^\prime(a-t+s) \,ds = R_\BX(a) -
R_\BX(a-t).
\]
In particular, if the process $\BX$ has a finite second spectral
moment, then the second derivative of the covariance function exists,
is continuous and negative at zero (unless the covariance function is
constant). Therefore, the derivative stays negative on an interval
around the origin, hence, the covariance function is concave on
$[0,a]$, and \eqref{econd1} holds, for $a>0$ small enough.

On the other hand, apart from degenerate cases, the situation
described in Proposition \ref{prcase1} cannot continue to hold for
arbitrarily large $a$. For example, if the covariance function
vanishes at infinity, then \eqref{econd1} fails for $a$ large enough and
$t=a/2$, say.

In addition, a simple calculation shows that it is always true that
%
\begin{eqnarray}
\label{simplecalcequn} \lim_{u\to\infty} \frac{1}{u^2} \log P \bigl(
\BX(0)>u, \BX(a)> u \bigr) =- \bigl(R_\BX(0)+R_\BX(a)
\bigr)^{-1}.
\end{eqnarray}
Combining this with \eqref{efirstC} shows that, in the scenario of
Proposition
\ref{prcase1}, the probability that $\BX$ exceeds a high level over
an entire interval
and the probability that it does so only at the endpoints of the
interval are, at a
logarithmic scale, the same.
\end{remark}

\begin{figure}

\includegraphics{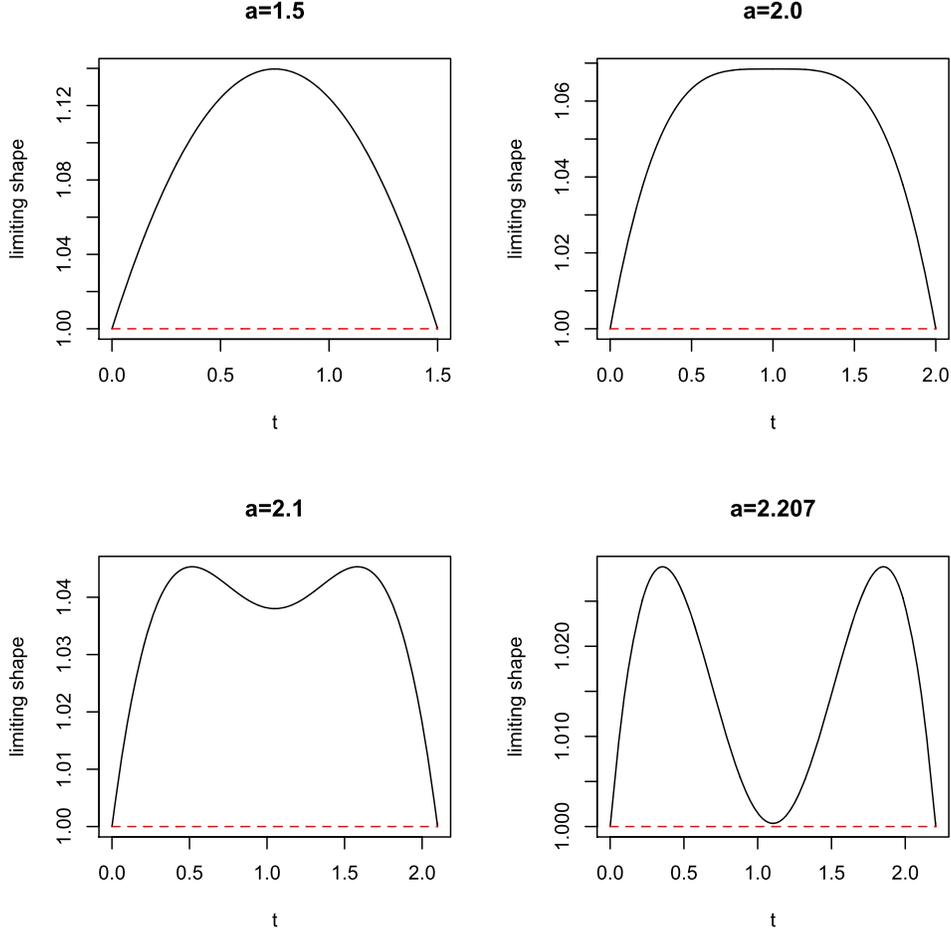}

\caption{Limiting shapes $x_a$ for the
stationary Gaussian process with covariance function $R_\BX(t)=\exp
(-t^2/2)$ when Proposition
\protect\ref{prcase1} applies.}\label{figcase1}
\end{figure}

The plots of Figure \ref{figcase1} show the limiting shape $x_a$ for the
stationary Gaussian process with covariance function
$R_\BX(t)=\exp(-t^2/2)$, for a range of $a$ for which Proposition
\ref{prcase1} applies. In this case the largest such $a$ is
approximately equal to~$2.2079$. See Example \ref{exGausscov} for
more details.

The plots of Figure \ref{figcase1} indicate that as $a$ approaches a
critical value (approximately $2.2079$ in this case), the limiting
curve $x_a$ ``attempts''
to cross
the level 1 at the midpoint of $[0,a]$. Equivalently, the normalized
process $u^{-1}\BX$
attempts to drop below level 1 at that point and so, speaking\vadjust{\goodbreak}
heuristically, it has to be
``supported'' at the midpoint $t=a/2$. The interpretation of Theorem
\ref{tfixedpath} in Remark \ref{rkinterpret} calls for adding a
mass to the measure $\mu^{(1)}$ for the critical value of $a$ at the
midpoint of the interval. The next result shows that, in certain
cases, this is indeed the optimal thing to do.

\begin{proposition} \label{prcase2}
Let $\BX$ be a stationary continuous Gaussian process. Suppose that,
for some $a>0$,
%
\begin{equation}
\label{emiddletoomuch} R_\BX(0)+R_\BX(a)>2
R_\BX(a/2) ,
\end{equation}
and let
%
\begin{equation}
\label{eepsilona} \vep_a=\frac{R_\BX(0)+R_\BX(a)-2 R_\BX(a/2)}{3R_\BX(0)+R_\BX(a)-4
R_\BX(a/2)} \in (0,1] .
\end{equation}
Suppose that for all $0\leq t\leq a/2$,
%
\begin{eqnarray}
\label{econd2} &&R_\BX(t)+R_\BX(a-t)- R_\BX(0)-R_\BX(a)
\nonumber
\\
&&\qquad\geq\vep_a \bigl[ R_\BX(t)+R_\BX(a-t)-
R_\BX(0)-R_\BX(a)\\
&&\hspace*{78pt}{} -2 \bigl( R_\BX(a/2-t) -
R_\BX(a/2) \bigr) \bigr] .\nonumber
\end{eqnarray}
Then a measure in $\Wa$ is given by
%
\begin{equation}
\label{esecondmu} \mu^{(2)} \definedas\frac{1-\vep_a}{2}
\delta_0+\frac{1-\vep_a}{2} \delta_1 +
\vep_a \delta_{1/2} .
\end{equation}
Furthermore,
%
\begin{eqnarray}
\label{esecondC} C_\BX(a)&= &\frac{3R_\BX(0)+R_\BX(a)-4
R_\BX(a/2)}{R_\BX(0)^2+R_\BX(0) R_\BX(a)-2 R_\BX(a/2)^2},
\\
\label{esecondh} h_a(x)& =& C_\BX(a) \biggl[
\frac{1-\vep_a}{2} \bigl( 1+e^{iax} \bigr) + \vep_a
e^{iax/2} \biggr]
\end{eqnarray}
for $F_\BX$-almost all $-\infty<x<\infty$, and
%
\begin{equation}
\label{esecondx} x_a(t) = C_\BX(a) \biggl[
\frac{1-\vep_a}{2} \bigl( R_\BX(t) + R_\BX(a-t) \bigr) +
\vep_aR_\BX\bigl(|t-a/2|\bigr) \biggr],
\end{equation}
$0\leq t\leq a$.
\end{proposition}
\begin{pf}
The proof is identical to that of Proposition \ref{prcase1} once we
observe that, under \eqref{emiddletoomuch}, $\mu^{(2)}$
is a legitimate probability measure.
\end{pf}

\begin{figure}

\includegraphics{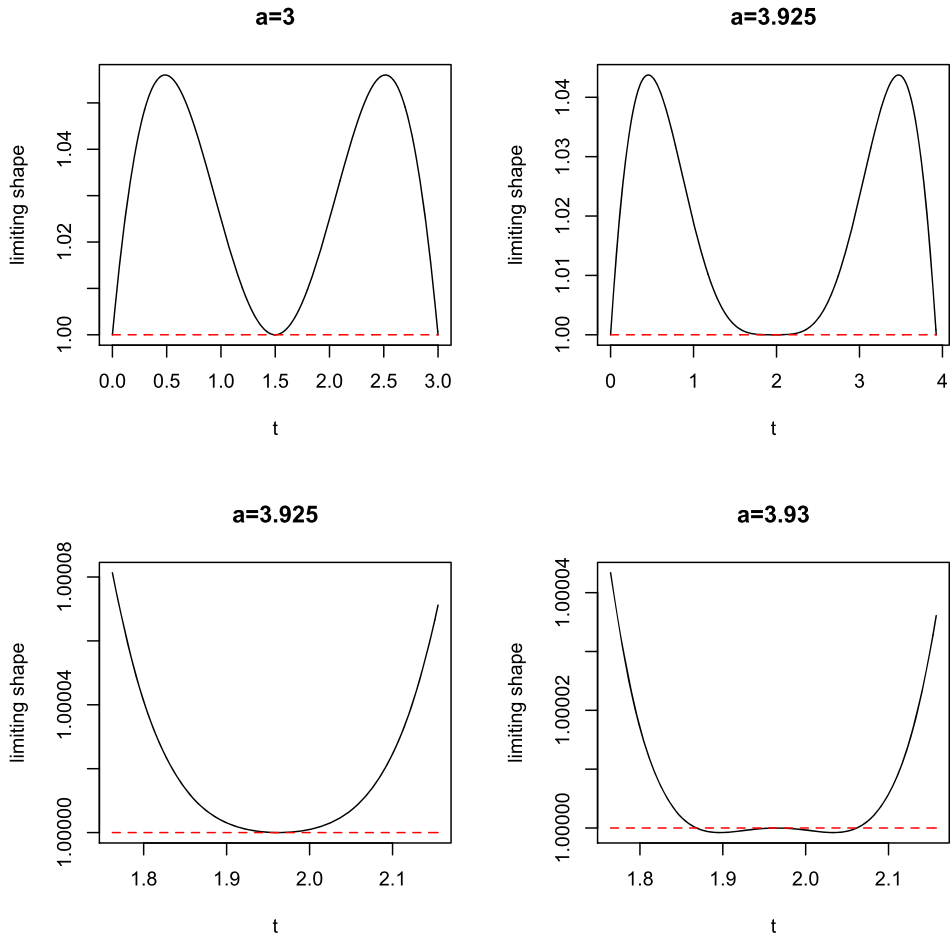}

\caption{Limiting shapes for $R_\BX(t)=\exp(-t^2/2)$ when Proposition
\protect\ref{prcase2} applies (the top row). The left plot in the
bottom row
is a blowup of the right plot in the top row. The right plot in the
bottom row shows how the constraints are violated soon after the
upper critical value of $a$.}
\label{figcase2}
\end{figure}

The plots of Figure \ref{figcase2} show the limiting shape $x_a$ for the
stationary Gaussian process with covariance function
$R_\BX(t)=\exp(-t^2/2)$, for a range of $a$ for which Proposition
\ref{prcase2} applies. In this case the range of $a$ is,
approximately, between $2.2079$ and $3.9283$. See Example
\ref{exGausscov} for more details.

\section{Specific covariance functions}
\label{examplessec}

In the previous section we saw some general results for one-dimensional
processes, with some
illustrative figures for what happens in the case of a Gaussian
covariance function.
In this section we look more carefully at this case, and also look at
what can be said for an
exponential covariance.

\begin{figure}

\includegraphics{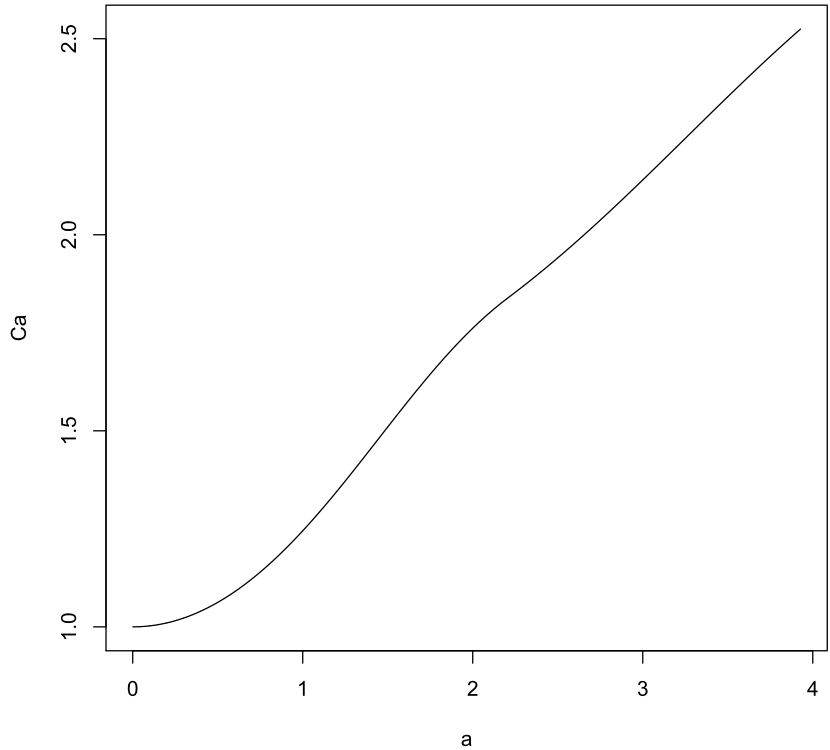}

\caption{The exponent $C_\BX(a)$ as a function of $a$ for $R_\BX
(t)=\exp(-t^2/2)$.}
\label{figCacase1}
\end{figure}

\begin{example} \label{exGausscov}
Consider the centered stationary Gaussian process with the Gaussian
covariance function
%
\begin{equation}
\label{eGausscov} R(t) = e^{-t^2/2},\qquad t\in\bbr.
\end{equation}
For this process the spectral measure has
a Gaussian spectral density which is of full support in
$\bbr$. In particular, for every $a>0$ there is a unique (symmetric)
measure of minimal energy. Furthermore, the second spectral moment is
finite, so that, according to Remark \ref{rkcase1}, for $a>0$
sufficiently small this process satisfies the conditions of
Proposition \ref{prcase1}. To find the range of $a$ for which this
happens, note that conditions \eqref{econd1} become, in this case,
%
\begin{equation}
\label{econdGauss1} e^{-t^2/2} + e^{-(a-t)^2/2}\geq1 + e^{-a^2/2},\qquad
0\leq t\leq a .
\end{equation}
Since the function
\[
g(t) = e^{-t^2/2} + e^{-(a-t)^2/2},\qquad 0\leq t\leq a ,
\]
is concave if $0\leq a\leq2$, and has a unique local minimum, at
$t=a/2$, when $a>2$, it is only necessary to check
\eqref{econdGauss1} at the midpoint $t=a/2$. At that point the
condition becomes
\[
\psi(a) = 2e^{-a^2/8}-1-e^{-a^2/2}\geq0 .
\]
The function $\psi$ crosses 0 at $a_1\approx2.2079$, which is the
limit of the validity of the situation of Proposition \ref{prcase1}
in this case. The plots of Figure \ref{figcase1} show the limiting
shape $x_a$ for this process in the situation of Proposition
\ref{prcase1}.

Somewhat longer (and numerical) calculations show that the conditions
of Proposition \ref{prcase2} hold for the process with the
covariance function \eqref{eGausscov} for an interval of values of
$a$ after the conditions of Proposition \ref{prcase1} break
down. The conditions of Proposition \ref{prcase2} continue to hold
until the second derivative at the midpoint $t=a/2$ of the limiting
function in
\eqref{esecondx} becomes negative (so that the function takes values
smaller than 1 in a neighborhood of the midpoint). To find when this
happens, we solve the equation
\[
\frac{1-\vep_a}{2} \bigl( R_\BX^{\prime\prime}(t) +
R_\BX^{\prime\prime}(a-t) \bigr) + \vep_aR_\BX^{\prime\prime}\bigl(|t-a/2|\bigr)=0
\]
at $t=a/2$. The resulting equation
\[
(1-\vep_a) \bigl(a^2/4-1 \bigr)e^{-a^2/8} -
\vep_a=0
\]
has the solution $a_2\approx3.9283$, which is the
limit of the validity of the situation of Proposition \ref{prcase2}
in this case. The plots of Figure \ref{figcase2} shed some light
on the above discussion. This discussion indicates, and calculations
confirm, that, in the next regime, the mass in the middle for the
optimal measure splits into two parts that start to move away from the
center. Heuristically, this is needed ``to support'' the trajectory
that, otherwise, would ``dip'' below 1 outside of the midpoint.

These
calculations rapidly become complicated. They seem to indicate that
the next regime continues to hold until around $a_3\approx5.4508$. In
this regime the optimal measure takes the form
%
\begin{equation}
\label{ethirdmu} \mu^{(3)} \definedas\frac{1-\vep_a}{2}
\delta_0+\frac{1-\vep_a}{2} \delta_1 + \frac{\vep_a}{2}
\delta_{1/2-d_a}+\frac{\vep_a}{2} \delta_{1/2+d_a},
\end{equation}
where $d_a$ is the distance of two internal masses from the midpoint.
When $a=4.5$, $\vep_{4.5}=0.36632$ and
$d_{4.5}=0.12285$, so that the internal atoms are at $0.37715$ and
$0.62285$, and the rest of the support is concentrated at the endpoints
of the interval with probabilities $0.31684$.
Figure \ref{figcase3} shows the limiting shape $x_{4.5}$.

\begin{figure}

\includegraphics{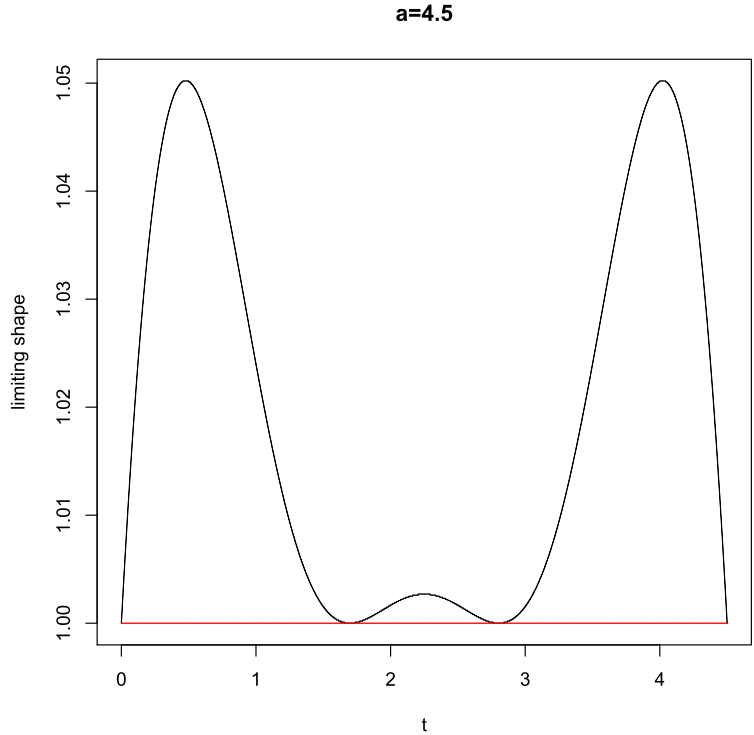}

\caption{The limiting shape in the case $a=4.5$ for $R_\BX(t)=\exp
(-t^2/2).$ }
\label{figcase3}
\end{figure}

It would be nice to understand all regimes, but we do not yet
know how to find a general structure. On the other hand, Section \ref{longsec}
gives asymptotic results for $a\to\infty$.

Finally, Figure \ref{figCacase1} shows the growth of the exponent
$C_\BX(a)$ with $a$ for as long as either Proposition \ref{prcase1}
or Proposition \ref{prcase2} applies.

\end{example}

The next example shows a situation very different from that of
Example~\ref{exGausscov}.\vadjust{\goodbreak}
%
\begin{example} \label{exO-U}
Consider an Ornstein--Uhlenbeck process, that is, a centered stationary
Gaussian process with the covariance function
%
\begin{equation}
\label{eO-Ucov} R(t) = e^{-|t|},\qquad  t\in\bbr.
\end{equation}
For this process the spectral measure has a Cauchy spectral density,
so it is also of full support in
$\bbr$. Therefore, for every $a>0$ there is a unique (symmetric)
measure of minimal energy. In this case, however, even the first
spectral moment is infinite. The covariance function is actually
\textit{convex} on the positive half-line so, in particular, the
conditions of
Proposition \ref{prcase1} fail for all $a>0$. In fact, it is
elementary to check that for the probability measure
%
\begin{equation}
\label{econtlebesgue} \mu=\frac{1}{a+2}\delta_0+
\frac{1}{a+2}\delta_1 + \frac
{a}{a+2}\lambda
\end{equation}
[where $\lambda$ is the Lebesgue measure on $(0,1)$], the integrals
\[
\int_0^1 R_\BX \bigl( a(u-v)
\bigr) \mu(du), \qquad 0\leq v\leq1 ,
\]
have a constant value, equal to $2/(a+2)$. Therefore, the measure
$\mu$ in \eqref{econtlebesgue} is the measure of minimal
energy, and $C_\BX(a)=(a+2)/2$ for all $a>0$.

By Theorem \ref{t1dimLDP} we conclude that the limiting function
$x_a$ is equal to 1 almost everywhere in $[0,a]$ with respect to the
Lebesgue measure. Since $x_a$ is continuous, it is identically equal
to 1 on $[0,a]$.
\end{example}

Examples \ref{exGausscov} and \ref{exO-U} demonstrate a number
of the ways a stationary Gaussian process ``prefers,'' in the large
deviations sense,
to stay above a high level
over an interval. The process of Example \ref{exGausscov} with
covariance function \eqref{eGausscov} is smooth; the most likely way
for it to stay above a level is to force it to be ``slightly'' above that
level at a properly chosen finite set of time points; after that it is
``held'' above the level at the rest of the interval $[0,a]$ by the
correlations of the process. The optimal configuration of the finite
set of points depends on the length of the interval $[0,a]$, and it
appears to undergo phase transitions at certain critical interval
lengths. The complete picture of this ``dynamical system'' of finite
sets remains unclear. On the other hand, the Ornstein--Uhlenbeck
process of Example \ref{exO-U} is continuous, but not smooth. In
fact, it behaves locally like a Brownian motion. Therefore,
``holding'' it ``slightly'' above a level at a discrete point does not
help, since it ``wants'' immediately to go below that
level. This explains the nature of the optimal measure $\mu$ in~\eqref{econtlebesgue},
and this nature stays the same no matter how
short or long the interval $[0,a]$ is. In particular, phase transitions
do not
happen for this process.

It remains to be investigated whether other types of behavior are
possible, and under what exact conditions on the Gaussian process each
type of behavior occurs. It is also likely that minimal energy
measures in $\Wa$ carry additional information, describing how
``slightly'' above the level $u$ a Gaussian process is most likely to
be, given that it is above that level along the interval. The exact
nature of this information also remains to be investigated.

\section{Asymptotics for long intervals}
\label{longsec}

In this section we investigate the asymptotics of the exponent
$C_\BX(a)$ for large $a$. We start with a result showing
that, for certain short memory stationary
Gaussian processes, the exponent $C_\BX(a)$ grows linearly with $a$
over long intervals. Furthermore, the energy of the uniform
distribution $\lambda$ on $[0,1]$ becomes, asymptotically,
minimal.\looseness=-1
%
\begin{theorem} \label{tLargeAShortMemory}
Let $\BX$ be a stationary continuous Gaussian process. Assume that
$R_\BX$ is positive, and satisfies the following condition:
%
\begin{equation}
\label{eshortmemory} \int_0^\infty R(t) \,dt<
\infty.
\end{equation}
Then, with $\lambda$ denoting the uniform probability measure on $[0,1]$,
%
\begin{eqnarray}
\label{elinearCx} \lim_{a\to\infty} \frac{1}{a}C_\BX(a) &=&
\biggl( \lim_{a\to\infty} a\int_0^1\int
_0^1 R_\BX \bigl( a(u-v) \bigr)
\lambda(du) \lambda(dv) \biggr)^{-1}
\nonumber
\\[-8pt]
\\[-8pt]
\nonumber
&=& \frac{1}{2\int_0^\infty R(t) \,dt} .
\end{eqnarray}
\end{theorem}
\begin{pf}
By Theorem \ref{t1dimLDP}, the statement of the present theorem
is equivalent to the following pair of claims:
%
\begin{equation}
\label{euniformlimit} \lim_{a\to\infty} a\int_0^1
\int_0^1 R_\BX \bigl( a(u-v)
\bigr) \lambda(du) \lambda(dv) = 2\int_0^\infty
R(t) \,dt
\end{equation}
and
%
\begin{equation}
\label{einfimumlimit}\qquad \liminf_{a\to\infty} a \min_{\mu\in
M_1^+([0,1])} \int
_0^1\int_0^1
R_\BX \bigl( a(u-v) \bigr) \mu(du) \mu(dv)\geq2\int
_0^\infty R(t) \,dt .
\end{equation}
Since
\begin{eqnarray*}
&&\int_0^1\int_0^1
R_\BX \bigl( a(u-v) \bigr) \lambda(du) \lambda(dv)
\\
&&\qquad = \frac{1}{a}\int_0^1 \biggl[ \int
_0^{av} R(t) \,dt + \int_0^{a(1-v)}
R(t) \,dt \biggr] \,dv ,
\end{eqnarray*}
\eqref{euniformlimit} immediately follows from
\eqref{eshortmemory} and the bounded convergence theorem. Therefore, it
only remains to prove \eqref{einfimumlimit}. Suppose that, to the
contrary, \eqref{einfimumlimit} fails, and choose a sequence
$a_n\to\infty$ such that
\[
\lim_{n\to\infty} a_n \min_{\mu\in
M_1^+([0,1])} \int
_0^1\int_0^1
R_\BX \bigl( a_n(u-v) \bigr) \mu(du) \mu(dv)< 2\int
_0^\infty R(t) \,dt .
\]
For each $n$ choose a symmetric $\mu_{a_n}\in{\mathcal W}_{a_n}$, so that
%
\begin{equation}
\label{ebada} \lim_{n\to\infty} a_n \int_0^1
\int_0^1 R_\BX \bigl(
a_n(u-v) \bigr) \mu_{a_n}(du) \mu_{a_n}(dv)< 2\int
_0^\infty R(t) \,dt .
\end{equation}
We claim that, for every $\gamma>0$,
%
\begin{equation}
\label{elowertail} \lim_{n\to\infty} \mu_{a_n} \bigl( \bigl[0,\gamma
a_n^{-1} \bigr] \bigr)=0 .
\end{equation}
Indeed, by the positivity of $R_\BX$, for any $\gamma>0$,
\begin{eqnarray*}
&&\int_0^1\int_0^1
R_\BX \bigl( a_n(u-v) \bigr) \mu_{a_n}(du)
\mu_{a_n}(dv)
\\
&& \qquad\geq \int_0^{\gamma a_n^{-1}}\int_0^{\gamma a_n^{-1}}
R_\BX \bigl( a_n(u-v) \bigr) \mu_{a_n}(du)
\mu_{a_n}(dv)
\\
&&\qquad \geq \bigl\{\mu_{a_n} \bigl( \bigl[0,\gamma a_n^{-1}
\bigr] \bigr) \bigr\}^2 \inf_{0\leq t\leq\gamma}R(t) ,
\end{eqnarray*}
so that \eqref{ebada} necessitates \eqref{elowertail}. Next,
define a sequence of signed measures on $[0,1]$ by
$\hat\mu_n = \mu_{a_n}-\lambda$. Note that
%
\begin{equation}
\label{ezeromass} \hat\mu_n \bigl( [0,1] \bigr) = 0\qquad \mbox{for each
$n$.}
\end{equation}
By the nonnegative definiteness of $R_\BX$,
\begin{eqnarray*}
&&\hspace*{-4pt}\int_0^1\int_0^1
R_\BX \bigl( a_n(u-v) \bigr) \mu_{a_n}(du)
\mu_{a_n}(dv)
\\
&&\hspace*{-6pt}\qquad = \int_0^1\int_0^1
R_\BX \bigl( a_n(u-v) \bigr) \lambda(du) \lambda(dv) + \int
_0^1\int_0^1
R_\BX \bigl( a_n(u-v) \bigr) \hat\mu_n(du)
\hat\mu_n(dv)
\\
&&\hspace*{-6pt}\qquad\quad{} + 2 \int_0^1\int_0^1
R_\BX \bigl( a_n(u-v) \bigr) \lambda(du) \hat
\mu_n(dv)
\\
&&\hspace*{-6pt}\qquad\geq\int_0^1\int_0^1
R_\BX \bigl( a_n(u-v) \bigr) \lambda(du) \lambda(dv) \\
&&\hspace*{-6pt}\qquad\quad{}+ 2 \int
_0^1\int_0^1
R_\BX \bigl( a_n(u-v) \bigr) \lambda(du) \hat
\mu_n(dv).
\end{eqnarray*}
We will show that
%
\begin{equation}
\label{enonneglimit} \lim_{n\to\infty} a_n \int
_0^1\int_0^1
R_\BX \bigl( a_n(u-v) \bigr) \lambda(du) \hat
\mu_n(dv)= 0 .
\end{equation}
Together with \eqref{euniformlimit} this will provide the necessary
contradiction to \eqref{ebada}. Let \mbox{$\gamma>0$}. Write the integral
in \eqref{enonneglimit} as
\begin{eqnarray*}
&&\int_{\gamma a_n^{-1}}^{1-\gamma a_n^{-1}} \biggl[ \int
_0^1 R_\BX \bigl(
a_n(u-v) \bigr) \,du \biggr] \hat\mu_n(dv)
\\
&&\quad{} + 2 \int_0^{\gamma a_n^{-1}} \biggl[ \int
_0^1 R_\BX \bigl(
a_n(u-v) \bigr) \,du \biggr] \hat\mu_n(dv)
\\
&&\qquad \definedas J_n^{(1)} + 2J_n^{(2)} .
\end{eqnarray*}
Observe that
\begin{eqnarray*}
\bigl| J_n^{(2)} \bigr| &=& \frac{1}{a_n} \biggl\llvert \int
_0^{\gamma a_n^{-1}} \biggl[ \int_0^{a_n v}
R_\BX(t) \,dt + \int_0^{a_n (1-v)}
R_\BX(t) \,dt \biggr] \hat\mu_n(dv) \biggr\rrvert
\\
&\leq& \frac{2 \int_0^{\infty} R_\BX(t) \,dt}{a_n} \|\hat\mu_n\| \bigl( \bigl[0,\gamma
a_n^{-1} \bigr] \bigr),
\end{eqnarray*}
so that by \eqref{elowertail} we obtain
%
\begin{equation}
\label{esmallJ2} \lim_{n\to\infty} a_n J_n^{(2)}
= 0
\end{equation}
for every $\gamma>0$. Next, we write
\begin{eqnarray*}
J_n^{(1)}&=& \frac{1}{a_n} \int_{\gamma a_n^{-1}}^{1-\gamma a_n^{-1}}
\biggl[ \int_0^{a_n v} R_\BX(t) \,dt +
\int_0^{a_n (1-v)} R_\BX(t) \,dt \biggr]
\hat \mu_n(dv)
\\
&=& \frac{2 \int_0^{\infty} R_\BX(t) \,dt}{a_n} \hat\mu_n \bigl( \bigl[\gamma
a_n^{-1}, 1-\gamma a_n^{-1} \bigr]
\bigr)
\\
&& {}- \frac{1}{a_n} \int_{\gamma a_n^{-1}}^{1-\gamma a_n^{-1}} \biggl[
\int_{a_n v}^\infty R_\BX(t) \,dt + \int
_{a_n (1-v)}^\infty R_\BX(t) \,dt \biggr] \hat
\mu_n(dv)
\\
&\definedas& J_n^{(11)}-J_n^{(12)}.
\end{eqnarray*}
It follows from \eqref{ezeromass} that
%
\begin{equation}
\label{esmallJ11}\qquad\bigl | a_n J_n^{(11)} \bigr| = 2 \int
_0^{\infty} R_\BX(t) \,dt\bigl | \hat
\mu_n \bigl( \bigl[0,\gamma a_n^{-1} \bigr]
\bigr) + \hat\mu_n \bigl( \bigl[1-\gamma a_n^{-1},1
\bigr] \bigr) \bigr|\to0
\end{equation}
as $n\to\infty$, by \eqref{elowertail}. Finally,
%
\begin{equation}
\label{esmallJ12} \bigl| a_n J_n^{(12)} \bigr| \leq4 \int
_\gamma^{\infty} R_\BX(t) \,dt,
\end{equation}
and we obtain by \eqref{esmallJ2}, \eqref{esmallJ11} and
\eqref{esmallJ12} that
\[
\limsup_{n\to\infty} a_n \biggl\llvert \int_0^1
\int_0^1 R_\BX \bigl(
a_n(u-v) \bigr) \lambda(du) \hat\mu_n(dv) \biggr\rrvert
\leq4 \int_\gamma^{\infty} R_\BX(t) \,dt.
\]
Letting $\gamma\to\infty$ proves \eqref{enonneglimit} and, hence,
completes the proof of the theorem.
\end{pf}

The next theorem is the counterpart of Theorem
\ref{tLargeAShortMemory} for certain long memory stationary Gaussian
processes. In this case,
the uniform distribution on $[0,1]$ is no longer, asymptotically,
optimal. We will assume that the covariance function of the process is
regularly varying at infinity:
%
\begin{equation}
\label{elongmemory} R_\BX(t)=\frac{L(t)}{|t|^\beta}, \qquad 0<\beta<1,
\end{equation}
where $L$ is slowly varying at infinity. Before stating the theorem,
we introduce new notation.

Consider the minimization problem
%
\begin{equation}
\label{eriesz} \min_{\mu\in M_1^+([0,1])} \int_0^1
\int_0^1 \frac{\mu(du)
\l\mu(dv)}{|u-v|^\beta},\qquad  0<\beta<1.
\end{equation}
This is a minimization problem of the same nature as in
\eqref{e1dimreprdual} with $a=0,b=1$, and the covariance function
$R_\BX$ replaced by the Riesz kernel $K_\beta(u,v)=|u-v|^{-\beta}$,
$u,v\in
[0,1]$. The general theory of energy of measures in
\cite{fuglede1960} applies to the Riesz kernel. In particular, the
minimum in \eqref{eriesz} is well defined, is finite and
positive. Let ${\mathcal W}_\beta$ be the set of measures in
$M_1^+([0,1])$ of minimal energy with respect to the Riesz
kernel. Note that the uniform measure $\lambda\notin{\mathcal
W}_\beta$ since it does not satisfy the optimality conditions in
Theorem 2.4 in~\cite{fuglede1960}.\looseness=-1

\begin{theorem} \label{tLargeALongMemory}
Let $\BX$ be a continuous stationary Gaussian
process. Assume that $R_\BX$ is positive and satisfies assumption
\eqref{elongmemory} of regular variation. Then for any
$\mu_\beta\in{\mathcal W}_\beta$,
%
\begin{eqnarray}
\label{enonlinearCx} \lim_{a\to\infty} R_\BX(a)\CXa = \biggl( \int
_0^1\int_0^1
\frac{\mu_\beta(du)
\l\mu_\beta(dv)}{|u-v|^\beta} \biggr)^{-1} .
\end{eqnarray}
\end{theorem}
\begin{pf}
Suppose first that there is a sequence $a_n\uparrow\infty$ such that
%
\begin{eqnarray}
\label{elowlim} &&\lim_{n\to\infty} \frac{1}{R_\BX(a_n)}\min_{\mu\in
M_1^+([0,1])} \int
_0^1\int_0^1
R_\BX \bigl( a_n(u-v) \bigr) \mu(du) \mu(dv)
\nonumber
\\[-8pt]
\\[-8pt]
\nonumber
&& \qquad <\int_0^1\int_0^1
\frac{\mu_\beta(du)
\mu_\beta(dv)}{|u-v|^\beta}.
\nonumber
\end{eqnarray}
For each $n$ choose $\mu_n\in{\mathcal W}_{a_n}$, let
$n_k\uparrow\infty$ be a subsequence such that $\mu_{n_k}\Rightarrow
\hat\mu$ weakly as $k\to\infty$ for some $\hat\mu\in
M_1^+([0,1])$. By Fatou's lemma and the regular variation of $R_\BX$,
\begin{eqnarray*}
&&\liminf_{k\to\infty} \frac{1}{R_\BX(a_{n_k})} \int_0^1
\int_0^1 R_\BX \bigl(
a_{n_k}(u-v) \bigr) \mu_{n_k}(du) \mu_{n_k}(dv)
\\
&&\qquad \geq\int_0^1\int_0^1
\frac{\hat\mu(du)
\hat\mu(dv)}{|u-v|^\beta}
 \geq\int_0^1\int_0^1
\frac{\mu_\beta(du)
\mu_\beta(dv)}{|u-v|^\beta},
\end{eqnarray*}
since $\mu_\beta$ has the smallest energy with respect to the Riesz
kernel. This contradicts~\eqref{elowlim}, thus proving that
%
\begin{eqnarray}
\label{elower} &&\liminf_{a\to\infty} \frac{1}{R_\BX(a)}\min_{\mu\in
M_1^+([0,1])} \int
_0^1\int_0^1
R_\BX \bigl( a(u-v) \bigr) \mu(du) \mu(dv)
\nonumber
\\[-8pt]
\\[-8pt]
\nonumber
&&\qquad\geq\int_0^1\int_0^1
\frac{\mu_\beta(du)
\mu_\beta(dv)}{|u-v|^\beta}.
\end{eqnarray}
In order to finish the proof, we need to establish a matching upper
limit bound.

To this end, let $\theta>0$ be a small number. We define a probability
measure $\nu_\theta\in M_1^+([0,1])$ by convolving $\mu_\beta$ with
the uniform distribution on $[0,\theta]$ and rescaling the resulting
convolution back to the unit interval. More explicitly, if $X$ and $U$ are
independent random variables, whose laws are $\mu_\beta$ and
$\lambda$, respectively, then $\nu_\theta$ is the law of $(X+\theta
U)/(1+\theta)$. Note that
%
\begin{equation}
\label{edensity} \nu_\theta\ll\lambda \qquad\mbox{and}\qquad \frac{d
\nu_\theta}{d\lambda}\leq
\frac{1+\theta}{\theta} \qquad\mbox{a.e. on $[0,1]$.}
\end{equation}

Given $0<\vep<1-\beta$, by Potter's bounds (see, e.g., Proposition 0.8
in \cite{resnick1987}), there is $t_0>0$ sufficiently large to ensure
\[
\frac{R_\BX(tx)}{R_\BX(t)}> (1-\vep)x^{-\beta-\vep}
\]
for all $t\geq t_0$ and $x\geq1$. We have
\begin{eqnarray*}
&&\min_{\mu\in
M_1^+([0,1])} \int_0^1\int
_0^1 R_\BX \bigl( a(u-v) \bigr)
\mu(du) \mu(dv)
\\
&& \qquad\leq\int_0^1\int_0^1
R_\BX \bigl( a(u-v) \bigr) \nu_\theta(du)
\nu_\theta(dv)
\\
&&\qquad = \int_0^1\int_0^1
\one \bigl( |u-v|\leq t_0/a \bigr) R_\BX \bigl( a(u-v) \bigr)
\nu_\theta(du) \nu_\theta(dv)
\\
&&\qquad\quad{} + \int_0^1\int_0^1
\one \bigl( |u-v|> t_0/a \bigr) R_\BX \bigl( a(u-v) \bigr)
\nu_\theta(du) \nu_\theta(dv)
\\
&&\qquad \definedas I_1(a) + I_2(a) .
\end{eqnarray*}
By the definition of $t_0$,
\[
\one \bigl( |u-v|> t_0/a \bigr) \frac{R_\BX (
a(u-v) )}{R_\BX(a)}\leq\frac{1}{1+\vep}|u-v|^{-(\beta+\vep)}
,
\]
so that by the dominated convergence theorem we have
\[
\lim_{a\to\infty} \frac{1}{R_\BX(a)} I_2(a) = \int
_0^1\int_0^1
\frac{\nu_\theta(du)
\nu_\theta(dv)}{|u-v|^\beta}.
\]
On the other hand, by \eqref{edensity},
\[
I_1(a)\leq R_\BX(0)\frac{2t_0}{a}\frac{1+\theta}{\theta}
= o \bigl(R_\BX(a) \bigr)
\]
as $a\to\infty$. We conclude that
\begin{eqnarray*}
&&\limsup_{a\to\infty}\frac{1}{R_\BX(a)}\min_{\mu\in
M_1^+([0,1])} \int
_0^1\int_0^1
R_\BX \bigl( a(u-v) \bigr) \mu(du) \mu(dv)
\\
&& \qquad \leq\int_0^1\int_0^1
\frac{\nu_\theta(du)
\nu_\theta(dv)}{|u-v|^\beta}.
\end{eqnarray*}
Once we show that
%
\begin{equation}
\label{enotheta} \lim_{\theta\to0} \int_0^1
\int_0^1 \frac{\nu_\theta(du)
\nu_\theta(dv)}{|u-v|^\beta} = \int
_0^1\int_0^1
\frac{\mu_\beta(du)
\mu_\beta(dv)}{|u-v|^\beta} ,
\end{equation}
we will have established an upper bound matching \eqref{elower}. This will
complete the proof of the theorem. Recall that \eqref{enotheta} is
equivalent to
\[
\lim_{\theta\to0} E \bigl| X_1-X_2+\theta(U_1-U_2)
\bigr|^{-\beta} = E | X_1-X_2 |^{-\beta} ,
\]
where $X_1,X_2,U_1,U_2$ are independent random variables, $X_1$ and
$X_2$ with the law $\mu_\beta$, while $U_1$ and $U_2$ are uniformly
distributed on $[0,1]$. This, however, follows by the dominated
convergence theorem and the following fact, that can be checked by
elementary calculations: there is $r_\beta\in(0,\infty)$ such that
for any $0<b<1$ and $0<\theta<1$,
\[
E \bigl| b+\theta(U_1-U_2) \bigr|^{-\beta}\leq
r_\beta b^{-\beta} .
\]
\upqed\end{pf}
%
\begin{remark} \label{rklrdconstant}
It follows from Proposition A.3 in \cite{khoshnevisanxiaozhong2003}
that the energy of the measure $\mu_\beta$ with respect to the Riesz
kernel cannot be smaller than one half of the energy of the uniform
measure. Hence,
\[
\lim_{a\to\infty} R_\BX(a)\CXa\in \bigl( (1-\beta) (2-\beta)/2,
(1-\beta) (2-\beta) \bigr) .
\]
\end{remark}

\section{The multidimensional case}
\label{secmultcase}

Our understanding of the one-dimensional case described in the
previous three sections,
while incomplete, is nevertheless quite significant.
In contrast, there is much less we can say about the
multivariate problem of Section \ref{secLD}. The problem lies, in
part, in the nonconvexity of the feasible set in
\eqref{eLDgenlimit} which leads, in turn, to the ``max-min'' problem
in Theorem \ref{tlargeenergypath}.

The following proposition is a multivariate version of Proposition
\ref{econd1}. Note that stationarity of the random field is not
required.
%
\begin{proposition} \label{prcase1gen}
Let $\BX= (X(\bt), \bt\in T)$ be a continuous Gaussian
random field on a compact set $T\subset\bbr^d$, and suppose that
$\ba, \bb$ are in $T$. Suppose that there is a path $\xi_0$ in $T$
connecting $\ba$ and
$\bb$ such that
%
\begin{equation}
\label{econd1gen}\qquad R_\BX \bigl(\ba,\xi_0(u) \bigr)+
R_\BX \bigl(\xi_0(u),\bb \bigr) \geq\frac{R_\BX(\ba,\ba)+2R_\BX(\ba,\bb)
+ R_\BX(\bb,\bb)}{2}>0
\end{equation}
for all $0\leq u\leq1$. Then the supremum in \eqref{ecapGenlimit} is
achieved on the path $\xi_0$ and
%
\begin{equation}
\label{eC1gen} \CXab= \frac{4}{R_\BX(\ba,\ba)+2R_\BX(\ba,\bb)
+ R_\BX(\bb,\bb)} .
\end{equation}
\end{proposition}
%
\begin{remark} \label{rksameends}
Using $u=0$ and $u=1$ in \eqref{econd1gen} shows that conditions of
Proposition \ref{prcase1gen} cannot be satisfied unless
$R_\BX(\ba,\ba)=R_\BX(\bb,\bb)$. Correspondingly, we can restate
\eqref{eC1gen} as
\[
\CXab= \frac{2}{R_\BX(\ba,\ba)+R_\BX(\ba,\bb)} .
\]
Recall \eqref{simplecalcequn}, which shows that this implies the
logarithmic equivalence
of the probabilities of $\BX$ being above the level $u$ along a curve
or at its endpoints.
\end{remark}
\begin{pf*}{Proof of Proposition \ref{prcase1gen}}
Consider the fixed path $\xi_0$. The assumption~\eqref{econd1gen}
shows that the measure
\[
\mu_0=\tfrac12 \delta_a + \tfrac12 \delta_b
\]
satisfies conditions \eqref{ehighestenergy} and, hence, is in
${\mathcal W}_{\xi_0}$ by Theorem
\ref{toptimalmchararct}. Therefore,
\begin{eqnarray*}
&&\min_{\mu\in M_1^+([0,1])} \int_0^1\int
_0^1 R_\BX \bigl(
\xi_0(u),\xi_0(v) \bigr) \mu(du) \mu(dv)
\\
&&\qquad = \int_0^1\int_0^1
R_\BX \bigl( \xi_0(u),\xi_0(v) \bigr)
\mu_0(du) \mu_0(dv)
\\
&& \qquad = \frac{R_\BX(\ba,\ba)+2R_\BX(\ba,\bb)
+ R_\BX(\bb,\bb)}{4} .
\end{eqnarray*}
On the other hand, for any other path in $T$ connecting $\ba$ and
$\bb$,
\begin{eqnarray*}
&&\min_{\mu\in
M_1^+([0,1])} \int_0^1\int
_0^1 R_\BX \bigl( \xi(u),\xi(v)
\bigr) \mu(du) \mu(dv)
\\
&&\qquad \leq \int_0^1\int_0^1
R_\BX \bigl( \xi(u),\xi(v) \bigr) \mu_0(du)
\mu_0(dv)
\\
&&\qquad = \frac{R_\BX(\ba,\ba)+2R_\BX(\ba,\bb)
+ R_\BX(\bb,\bb)}{4} .
\end{eqnarray*}
Therefore, the supremum in \eqref{ecapGenlimit} is
achieved on the path $\xi_0$, and \eqref{eC1gen} follows
by Theorem \ref{tlargeenergypath}.
\end{pf*}

Even for the most common Gaussian random fields, the assumptions
of Proposition \ref{prcase1gen} may be satisfied on some path but
not on the straight line connecting the two points. In that case, the
straight line, clearly, fails to be optimal.

\begin{example} \label{exbrowniansheet}
Consider a Brownian sheet in $d\geq2$ dimensions. This is the
continuous centered Gaussian random field $\BX$ on $[0,\infty)^d$ with
covariance function
\[
R_\BX(\bs,\bt) = \prod_{j=1}^d
\min(s_j,t_j), \qquad\bs,\bt\in [0,\infty)^d .
\]
We restrict the random field to the hypercube $T=[0,d]^d$, and let
\[
\ba= (1,2,\ldots,d-1, d),\qquad \bb=(d,1,2,\ldots, d-1) .
\]
It is elementary to check that the path
\[
\xi_0(u) = \cases{ %
\bigl( 1+d(d-1)u,2,
\ldots, d-1,d \bigr), & \quad$ \mbox{for } 0\leq u\leq \displaystyle\frac{1}{d},$
\vspace*{2pt}\cr
( d,1,\ldots, j-2, 2j-1-du,j+1,\ldots, d ), & \quad $\mbox{for } \displaystyle\frac{j-1}{d}\leq u
\leq\frac{j}{d},$
\vspace*{2pt}\cr
& \quad $j=2,\ldots, d$, }
\]
satisfies \eqref{econd1gen} and, hence, the supremum in \eqref
{ecapGenlimit} is
achieved on that path. Therefore, by Proposition \ref{prcase1gen},
\[
\CXab= \frac{2}{d! + (d-1)!} .
\]
On the other hand, if we consider the straight line connecting the
points $\ba$ and $\bb$,
\[
\xi(u) = \bigl( 1+(d-1)u,2-u,3-u, \ldots, d-u \bigr),\qquad 0\leq u\leq 1 ,
\]
then the sum in the right-hand side of \eqref{econd1gen} becomes
\[
L(u) = \prod_{j=2}^d (j-u) + \bigl(
1+(d-1)u \bigr) (d-1)!,\qquad  0\leq u\leq1 .
\]
The function $L$ achieves the value $d!+(d-1)!$ at the endpoints $u=0$
and \mbox{$u=1$}, and is strictly convex if $d\geq3$. Therefore, it takes
values strictly smaller than $d!+(d-1)!$ over $0<u<1$. That
is, \eqref{econd1gen} fails, and the straight line is not
optimal. If $d=2$, however, then $L$ is a constant function, condition
\eqref{econd1gen} holds over the straight line path, and the
straight line is optimal.

We also note that, if $d=1$, then the Brownian sheet becomes the
Brownian motion in one dimension. In that case it is, clearly,
impossible to find two positive points $a<b$ in which the process has
the same variance, so Proposition \ref{prcase1gen} does not
apply. In this case, however, we are in the situation of Theorem
\ref{t1dimLDP}, so if $0<a<b<\infty$, then the measure
$\mu=\delta_a$ satisfies \eqref{echeck1dim} and, hence, is optimal.
\end{example}

The above example notwithstanding, under certain assumptions on the
random field, the straight line path between two points turns out to
be optimal for the optimization problem
\eqref{ecapGenlimit}. The next result describes one such
situation.

Recall that a random field on $\bbr^d$ is
\textit{isotropic} if its law is invariant under rigid motions of the
parameter space. A centered Gaussian random field $\BX$ is isotropic
if and only if its covariance function is a function of the Euclidian
distance between two points. With the usual abuse of notation we will
write $R_\BX(\ba,\bb)=R_\BX(\|\bb-\ba\|)$, $\ba,\bb\in\bbr^d$.
%
\begin{proposition} \label{prisotropic}
Let $\BX$ be a continuous centered isotropic Gaussian random field,
such that the covariance function $R_\BX$ is nonincreasing. Then for
any $\ba,\bb\in\bbr^d$, the straight path connecting the points
$\ba$
and $\bb$ is optimal for the optimization problem
\eqref{ecapGenlimit}.
\end{proposition}
\begin{pf}
We may and will assume, without loss of generality, that
$\ba=(a,0,\ldots,0)$ and $\bb=\bnull$ for some $a>0$.
We start with showing that the supremum over $\xi\in\cP(\bnull,\ba)$
is achieved over paths in
\[
\cP_l= \bigl\{ \xi\dvtx [0,1]\to \bigl\{ (x,0,\ldots, 0)\dvtx x
\geq0 \bigr\}, \mbox{continuous, } \xi(0)=\bnull, \xi(1) =\ba \bigr\} .
\]
To this end, it is enough to show that for each $\xi\in
\cP(\bnull,\ba)$ there is $\hat\xi\in\cP_l$ such that
%
\begin{eqnarray}
\label{edominatepath} &&\min_{\mu\in
M_1^+([0,1])} \int_0^1
\int_0^1 R_\BX \bigl( \bigl\|\xi(u)-
\xi(v)\bigr\| \bigr) \mu(du) \mu(dv)
\nonumber
\\[-8pt]
\\[-8pt]
\nonumber
&&\qquad  \leq\min_{\mu\in
M_1^+([0,1])} \int_0^1\int
_0^1 R_\BX \bigl( \bigl\|\hat\xi(u)-\hat
\xi(v)\bigr\| \bigr) \mu(du) \mu(dv) .
\end{eqnarray}
To see this, define for $\xi\in\cP(\bnull,\ba)$
\[
\hat\xi(u) = \bigl(\bigl \|\xi(u)\bigr\|,0,\ldots,0 \bigr),\qquad 0\leq u\leq1 .
\]
Clearly, $\hat\xi\in\cP_l$, and \eqref{edominatepath} follows by
the monotonicity of $R_\BX$ and the triangle inequality
\[
\bigl| \bigl\|\xi(u)\bigr\| - \bigl\|\xi(v)\bigr\| \bigr| \leq\bigl\|\xi(u)-\xi(v)\bigr\| .
\]

Next, any $\xi\in\cP_l$ is of the form
%
\begin{equation}
\label{e1dimpath} \xi(u) = \bigl(\varphi(u),0,\ldots,0 \bigr), \qquad 0\leq u\leq1 ,
\end{equation}
with $\varphi\dvtx  [0,1]\to[0,\infty)$ a continuous function,
satisfying $\varphi(0)=0$, $\varphi(1)=a$. Defining
$\hat\varphi(u)=\min(\varphi(u),a)$, $0\leq u\leq1$, and
\[
\hat\xi(u) = \bigl(\hat\varphi(u),0,\ldots,0 \bigr), \qquad 0\leq u\leq1 ,
\]
we see that the supremum over paths in $\cP_l$ is, actually, achieved
over paths whose image is exactly the interval $[0,a]$. Finally, for
any path $\xi\in\cP_l$ of the latter type, given in the form
\eqref{e1dimpath}, define
\[
r(v) = \inf \bigl\{ u\in[0,1]\dvtx \varphi(u)=av \bigr\},\qquad 0\leq v\leq 1 .
\]
Then $r$ is a measurable map from $[0,1]$ to itself, so for
any $\mu\in
M_1^+([0,1])$, we can define $\mu_1\in M_1^+([0,1])$ by $\mu_1=\mu
\circ
r^{-1}$. Then
\begin{eqnarray*}
&&\int_0^1\int_0^1
R_\BX \bigl( \bigl|\varphi(u)-\varphi(v)\bigr| \bigr) \mu_1(du)
\mu_1(dv)
\\
&&\qquad =\int_0^1\int_0^1
R_\BX \bigl( \bigl|\varphi \bigl(r(u) \bigr)-\varphi \bigl(r(v) \bigr)\bigr|
\bigr) \mu(du) \mu(dv)
\\
&&\qquad = \int_0^1\int_0^1
R_\BX\bigl ( a|u-v| \bigr) \mu(du) \mu(dv).
\end{eqnarray*}
Therefore,
\begin{eqnarray*}
&&\min_{\mu\in
M_1^+([0,1])} \int_0^1\int
_0^1 R_\BX \bigl( \bigl\|\xi(u)-\xi(v)\bigr\|
\bigr) \mu(du) \mu(dv)
\\
&&\qquad \leq \min_{\mu\in M_1^+([0,1])} \int_0^1\int
_0^1 R_\BX \bigl( a|u-v| \bigr) \mu(du)
\mu(dv),
\end{eqnarray*}
and the statement of the proposition follows.
\end{pf}

%
%

%
%

%



\printaddresses

\end{document}